\documentclass{article}
\usepackage{tikz}
\usepackage[left=3cm,right=3cm,top=4cm,bottom=3.5cm]{geometry}
\usetikzlibrary{arrows,shapes,automata,petri,positioning,calc}
\usepackage{amsmath}
\usepackage{amsthm}
\theoremstyle{definition}
\newtheorem{definition}{Definition}
\usepackage{graphicx}
\usepackage{wrapfig}
\usepackage{caption}
\usepackage{csquotes}
\usepackage{subcaption}
\usepackage{mathtools}
\usepackage{sidecap}
\usepackage{hyperref}
\usepackage[title]{appendix}%
\usepackage{algorithm}
\usepackage{bbm}
\usepackage{lmodern,bm}
\usepackage{amsfonts}
\usepackage{fancyhdr}
\usepackage{mathrsfs}
\usepackage{float}
\usepackage{enumerate}
\usepackage{enumitem}
\usepackage[english]{babel}
\usepackage[nottoc]{tocbibind}
\usepackage{authblk}
\usepackage{physics}

\usepackage[backend=biber,giveninits=true,bibstyle=authoryear,url = false,doi = false,isbn=false]{biblatex}
\makeatletter
\newcommand*{\diff}%
{\@ifnextchar^{\DIfF}{\DIfF^{}}}
\def\DIfF^#1{%
  \mathop{\mathrm{\mathstrut d}}%
\nolimits^{#1}\gobblespace}
\def\gobblespace{%
\futurelet\diffarg\opspace}
\def\opspace{%
  \let\DiffSpace\!%
  \ifx\diffarg(%
    \let\DiffSpace\relax
    \else
    \ifx\diffarg[%
      \let\DiffSpace\relax
      \else
      \ifx\diffarg\{%
        \let\DiffSpace\relax
      \fi\fi\fi\DiffSpace}

\DeclareNameAlias{sortname}{family-given}
\DeclareNameAlias{default}{family-given}

\DeclareFieldFormat[inbook]{title}{#1}
\DeclareFieldFormat[article]{title}{#1}
\DeclareFieldFormat[inproceedings]{title}{#1}

\defbibenvironment{bibliography}
  {\enumerate
     {}
     {\setlength{\labelwidth}{\labelnumberwidth}%
      \setlength{\leftmargin}{\labelwidth}%
      \setlength{\labelsep}{\biblabelsep}%
      \addtolength{\leftmargin}{\labelsep}%
      \setlength{\itemsep}{\bibitemsep}%
      \setlength{\parsep}{\bibparsep}}%
      }
  {\endlist}
  {\item}
\title{{Performance Evaluation of Small Call Centres in a Transient Regime}}
\author[1,2]{Mark Fackrell}
\author[1,2]{Hritika Gupta\thanks{\href{mailto:hritikag@student.unimelb.edu.au}{hgupta@uow.edu.au}}}
\author[1,2]{Peter G. Taylor}
\affil[1]{School of Mathematics and Statistics, The University of Melbourne}
\affil[2]{ARC Training Centre in Optimisation Technologies, Integrated Methodologies, and Applications (OPTIMA)}

\begin{document}

\maketitle

\begin{abstract}
    {This paper addresses a fundamental and practically significant problem in call centre operations -- determining optimal call allocation policies that meet client service targets while minimising staffing costs. Motivated by a problem presented by an industry partner, we examine a real-world setting involving a relatively small call centre with hierarchical structure among agents. It is natural to model the operation of such a centre as a continuous-time Markov chain.
    
    To gain insight into the structure of optimal policies, we first (i) apply backward induction based on Bellman’s equation to a finite-horizon discrete-time model, and (ii) derive stationary policies for an infinite-horizon continuous-time model with discounting.
    
    Subsequently, we evaluate the performance of these policies in the original finite-horizon continuous-time setting by computing the expected number of abandonments and the waiting time distributions of customers. This is achieved using first-step analysis combined with Laplace transform methods. The effectiveness of the proposed approach is illustrated through numerical examples. }
\end{abstract}

\section{Introduction}\label{intro}

The term \textit{call centre} is commonly used to describe a telephone-based human-service operation. A call centre provides teleservices in which customers and service agents are remote from each other (\textcite{mandelbaum2000empirical}). They have become increasingly popular, and service quality is now a crucial factor in attracting and retaining customers.

Call centre service providers frequently have contracts with their clients under which they have to manage their client's customer care department and meet certain quality standards, such as serving at least $90\%$ of the calls within 20 seconds or having less than $3\%$ of call abandonments per month. Human resources typically account for more than $60\%$ of the total operating cost of a call centre (\textcite{gans_telephone_2003}, \textcite{brown_statistical_2005}). Hence, a reasonable objective for a call centre service provider is to minimise the staffing requirements while meeting all the targets. This can be achieved by improving the prediction accuracy related to call arrival processes, using an efficient agent scheduling method and making the best use of resources by appropriately allocating incoming calls to agents. It is the third of these that we will focus on in this paper.


In many call centres, customer queries are divided into different categories depending on the level of complexity of the query. Similarly, agents are also assigned a skill level depending on their ability to solve these queries. When a call arrives, based on the caller's initial input to an Interacting Voice Response (IVR) system, they are assigned to an agent with the appropriate skill level. Some agents are multi-skilled and can deal with more than one type of query. When a call centre has multi-skilled agents, there can be a number of possibilities for allocating an incoming call to an agent. We call the rules of assigning an incoming call to an agent a \textit{policy}.

Skills-based routing has been studied many times using different techniques. 
\textcite{wallace_staffing_2005} developed a suboptimal routing algorithm for skills-based routing, and showed using simulation analysis that limited cross training among agents could be almost as efficient as training all the agents for all the skills.
\textcite{stolletz2004performance} calculated the steady state probabilities and performance measures for a given priority based routing policy. 
\textcite{bhulai2009dynamic} obtained nearly optimal routing policies by minimising the long term expected average cost for the specific case where there is either no buffer for waiting customers or there is a buffer and all agents have either a single skill or are fully cross-trained.
\textcite{mehrbod2018caller} and \textcite{jorge2020intelligent} used machine learning techniques such as Random Forests and Neural Networks for their analysis. Simulation, see for example \textcite{adetunji2007performance}, \textcite{akhtar2010exploiting} and \textcite{mehrotra2012routing}, has also been a popular method for comparing different routing policies.

\subsection{Motivation}

Our study is motivated by a problem brought to us by an industry partner that manages call centres for its clients. The company was struggling to meet targets in terms of expected numbers of abandonments and losses, despite having high non-occupancy rates among its agents. It wanted us to study a particular client's call centre that was reasonably small, with around a dozen agents in total.

While our industry partner utilised models to estimate its workload and determine the required total number of agents, it did not use mathematically-backed strategies for allocating calls to agents. Indeed, call allocation was managed by an employee who manually decided how to allocate calls based on observing the system.

{Table \ref{tab:callcentre-snapshot} presents a sample dataset that is representative of the operational data received from a multi-tier call centre on January 3, 2022. There are four levels of agents of increasing expertise who were required to handle four levels of queries of increasing complexity. The data are reported at thirty-minute intervals for four service levels. For each level, the table records the number of calls offered, answered, and abandoned, together with the proportion of abandoned calls, the service level (percentage of calls answered within the target threshold), the average speed of answer (ASA), and the average handling time (AHT). The data illustrate common operational characteristics of call centres, such as higher abandonment rates and lower service levels during periods of heavy demand, and longer handling times for higher service levels where queries tend to be more complex. The data also demonstrates that, while the arrival processes could reasonably be modelled as homogeneous over a half-hourly time-scale, they were inhomogeneous over longer time scales.}

\begin{table}[!h]
\centering
\resizebox{0.95\textwidth}{!}{%
\begin{tabular}{lccccccc}
\hline
\textbf{Interval} & \textbf{Offered} & \textbf{Answered} & \textbf{Abandoned} & \textbf{Abandon \%} & \textbf{Service Level \%} & \textbf{ASA} & \textbf{AHT} \\ 
\hline
\multicolumn{8}{l}{\textbf{Level 1}} \\
8:00 AM & 10 & 10 & 0 & 0\% & 100\% & 00:00:03 & 00:02:50 \\
8:30 AM & 20 & 19 & 1 & 5\% & 60\% & 00:00:26 & 00:03:05 \\
9:00 AM & 12 & 12 & 0 & 0\% & 92\% & 00:00:04 & 00:06:05 \\
9:30 AM & 10 & 10 & 0 & 0\% & 90\% & 00:00:03 & 00:03:10 \\
10:00 AM & 11 & 11 & 0 & 0\% & 90\% & 00:00:07 & 00:03:12 \\
10:30 AM & 15 & 13 & 2 & 13\% & 45\% & 00:00:48 & 00:04:55 \\
11:00 AM & 20 & 20 & 0 & 0\% & 80\% & 00:00:33 & 00:03:40 \\
11:30 AM & 12 & 12 & 0 & 0\% & 100\% & 00:00:02 & 00:05:20 \\
12:00 PM & 18 & 18 & 0 & 0\% & 94\% & 00:00:05 & 00:04:00 \\

\multicolumn{8}{l}{\textbf{Level 2}} \\
8:00 AM & 16 & 16 & 0 & 0\% & 100\% & 00:00:05 & 00:06:30 \\
8:30 AM & 25 & 24 & 1 & 4\% & 64\% & 00:00:29 & 00:06:10 \\
9:00 AM & 31 & 30 & 1 & 3\% & 88\% & 00:00:10 & 00:06:18 \\
9:30 AM & 40 & 37 & 3 & 7\% & 53\% & 00:01:05 & 00:07:05 \\
10:00 AM & 55 & 44 & 11 & 20\% & 0\% & 00:04:20 & 00:07:10 \\
10:30 AM & 50 & 43 & 7 & 14\% & 2\% & 00:04:10 & 00:06:25 \\
11:00 AM & 36 & 36 & 0 & 0\% & 41\% & 00:00:58 & 00:07:35 \\
11:30 AM & 53 & 49 & 0 & 0\% & 81\% & 00:00:21 & 00:05:55 \\
12:00 PM & 45 & 45 & 0 & 0\% & 61\% & 00:00:40 & 00:07:45 \\

\multicolumn{8}{l}{\textbf{Level 3}} \\
8:00 AM & 9 & 8 & 0 & 0\% & 75\% & 00:00:23 & 00:05:18 \\
8:30 AM & 14 & 13 & 1 & 7\% & 52\% & 00:00:48 & 00:03:30 \\
9:00 AM & 18 & 18 & 0 & 0\% & 94\% & 00:00:09 & 00:06:10 \\
9:30 AM & 22 & 20 & 2 & 9\% & 30\% & 00:01:42 & 00:04:20 \\
10:00 AM & 25 & 24 & 1 & 4\% & 19\% & 00:02:12 & 00:06:55 \\
10:30 AM & 22 & 22 & 0 & 0\% & 18\% & 00:02:05 & 00:06:25 \\
11:00 AM & 18 & 18 & 0 & 0\% & 44\% & 00:01:09 & 00:05:35 \\
11:30 AM & 23 & 22 & 1 & 4\% & 55\% & 00:01:05 & 00:05:50 \\
12:00 PM & 24 & 23 & 2 & 8\% & 35\% & 00:01:09 & 00:06:40 \\

\multicolumn{8}{l}{\textbf{Level 4}} \\
8:00 AM & 3 & 3 & 0 & 0\% & 50\% & 00:01:36 & 00:20:45 \\
8:30 AM & 3 & 3 & 0 & 0\% & 25\% & 00:01:09 & 00:13:50 \\
9:00 AM & 7 & 6 & 0 & 0\% & 75\% & 00:00:15 & 00:06:20 \\
9:30 AM & 5 & 5 & 1 & 20\% & 14\% & 00:02:52 & 00:18:40 \\
10:00 AM & 5 & 5 & 1 & 20\% & 0\% & 00:05:25 & 00:09:35 \\
10:30 AM & 2 & 2 & 0 & 0\% & 0\% & 00:06:25 & 00:19:40 \\
11:00 AM & 4 & 3 & 1 & 25\% & 25\% & 00:02:25 & 00:17:05 \\
11:30 AM & 3 & 3 & 0 & 0\% & 100\% & 00:00:02 & 00:12:15 \\
12:00 PM & 5 & 3 & 2 & 40\% & 0\% & 00:02:43 & 00:09:50 \\
\hline
\end{tabular}
}
\caption{Snapshot of call centre data on January 3, 2022}
\label{tab:callcentre-snapshot}
\end{table}

{There is an extensive literature on call centre performance in the Halfin-Whitt regime, which focuses on systems operating under many-server heavy-traffic conditions. In contrast, our goal is to develop performance measurement tools that do not rely on the assumption that the call centre is large. The analysis we propose is valid regardless of the system size, however, the computational component is particularly well suited to small call centres due to the curse of dimensionality that arises in larger systems.}

We thought that a natural analysis would be to use a sequence of finite-horizon continuous-time Markov decision models for each half-hour interval. Since it is unreasonable that the state at the end of one interval would follow the stationary distribution of the next, such an analysis would necessarily have to use transient distributions and measures.

However, transient analysis of finite-horizon continuous-time Markov decision processes is difficult to carry out because the number of transitions occurring within the time horizon is both random and policy-dependent (\textcite{van2018uniformization}).  

A widely used technique for the transient analysis of continuous-time Markov chains is uniformisation. For a continuous-time Markov chain with transition rate matrix $Q$ that has bounded diagonal entries, this technique essentially writes the matrix exponential in the form
\begin{equation}
    \label{eq:uniformisation}
e^{Qt} = e^{-\lambda t}\sum_{k=0}^\infty \frac{(\lambda t)^k}{k!} P^k
\end{equation}
where $\lambda \geq \sup_{i} |q_{ii}|$ and $P=Q/\lambda + I$ with $I$ the identity matrix (see \textcite{gross1984randomization}, \textcite{buchholz2011numerical}).

The matrix $P$ on the right hand side of Equation {\eqref{eq:uniformisation} is the transition matrix of a discrete-time Markov chain, which facilitates a discrete-time analysis.
Uniformisation can be incorporated into a Markov decision process by letting the entries of $P$ depend on state-dependent actions. While this approach can be effective for many problems, it is not the most efficient approach for our problem. Observe the Markov chain with transition matrix $P$ frequently has states with a high probability of `self-transition' back to the same state. Putting that together with the above observation that, in our setting, decisions are required only at certain transition points, applying uniformisation would introduce a large number of redundant decision points, significantly increasing the computational burden. Since our goal is to develop an optimal policy that is robust, scalable, and computationally practical for a real-world call centre with multi-level agent hierarchies, we decided not to employ uniformisation.

The approach that we discuss below is related to the framework presented in \textcite[Sections 10.5 and 10.6]{kulkarni2014modeling}. Kulkarni discussed semi-Markov decision processes (SMDPs), where decisions are made at transition epochs occurring at random times, and the system evolves according to continuous-time dynamics between transitions. Our model is Markovian, but this approach still works. 

We first analyse
\begin{enumerate}
\item the discrete-time jump chain embedded in the continuous time process over a fixed finite number of transitions and
\item an approximating continuous-time Markov decision process over an infinite horizon with discounting.
\end{enumerate}
For the first model, we make use of backward induction using Bellman's equation to characterise the value function and derive an optimal policy. In doing this, we have ignored the variability in the sojourn times that occurs in the original continuous-time model. {In doing this, we are also assuming that a decision can only be taken at a state transition.} \textcite{koole2015optimization} used a similar analysis for optimal policies that depend on the waiting time of the longest waiting call, using the steady-state distribution of the underlying continuous-time Markov chain to approximate the waiting time distributions of the customers. A key question that remains is how to match the bound on the number of transitions in the discrete-time model to a finite-horizon analysis in continuous-time. 

Our analysis of the second model is necessarily stationary. However, again we can expect to gain insight into the structure of policies that are optimal for the original model.

Based on the hierarchical structure of our industry partner's centre, we first consider a model with $\nu$ classes of query and $\nu$ levels of agents, where level $i$ agents have the skills to deal with level $i$ and level $i-1$ queries. We work with this model in Sections \ref{General}, \ref{ch:res}, and \ref{ab4}.

To illustrate the approach, in Section \ref{General} we take $\nu = 2$ and use Bellman's equation to find the optimal policy for the embedded discrete-time Markov decision process over a finite horizon and the discounted continuous-time Markov decision process over an infinite horizon. This gives us an indication of the class of policies that are optimal. We discuss this in Section \ref{ch:res}.

{In Section \ref{ab4}, we revert back to the model of our industry partner's centre that has $\nu = 4$. Adapting the method of \textcite{chiera_what_2002}, we then develop a method to calculate the expected number of abandonments and waiting time distributions in call centres that can be modelled as continuous-time Markov chains under a transient setting, without using any asymptotic regimes. 

The expected number of abandonments and the distribution function of customer waiting times are the two main quantities required to assess performance measures specified in service level agreements between call centre providers and their clients. Computing distribution functions is generally more challenging than computing expectations, particularly when there are multiple types of callers. Depending on the policy, a newly arriving customer may be allowed to jump ahead of others in the queue, meaning that a customer’s waiting time can be influenced by subsequent arrivals.

There has been substantial research on waiting time distributions for queueing systems under stationary conditions, but relatively few studies focus on their transient behaviour. Transient waiting time distributions have primarily been investigated for single-server queues only.
\textcite{blomqvist1970transient} studied the standard deviation and covariance function of transient waiting time distributions for the $GI/G/1$ queue. 
\textcite{dalen1980transient} derived an integral solution for the waiting time distributions for the $G/M/1$ priority queue.
\textcite{bertsimas1992transient} derived the closed form expressions for the Laplace transforms of the queue length distribution and the waiting time distribution for the $G/G/1$ queue.
\textcite{abate1993calculation} developed an algorithm to numerically invert multidimensional transforms and applied it to find the transient queue length and workload distribution of the $M/G/1$ queue.
Lower bounds for the tails of both the steady-state and transient waiting time distributions for the $M/GI/s$ queue were established in \textcite{whitt2000impact}. 
\textcite{van2005approximated} proposed a method to approximate transient performance measures of a discrete-time queueing system via steady-state analysis. 

\textcite{knessl2015transient} derived explicit expressions for the Laplace transform of the time-dependent distribution and the first passage time of the $M/M/m + M$ queue (an $M/M/m$ queue with abandonments in which customer patience times are exponentially distributed) using Green’s functions and contour integrals related to hypergeometric functions. Although these expressions provide valuable analytical insight, their complexity makes them difficult to implement computationally. This study is the closest to the queueing system we consider, namely the $M/M/k/\ell + M$ model. However, our work extends this framework by incorporating multiple levels of agents in a finite-capacity system and a reservation policy, which further complicates the Laplace transform expressions of our queues. Hence, we choose to numerically invert the Laplace transform expressions. 

Discussions with our call centre partner revealed that many call centres still rely on the Erlang-A formula in their operations, despite recognising that its assumptions are not well suited to real-world, complex systems. Although there is extensive literature on optimal allocation policies for call centres, most studies are either highly analytical and too technical for companies to implement, or entirely based on simulations and simplified queueing models. The aim of this research is to bridge that gap by producing results that do not rely on simulations, black-box systems, or random components, but instead yield exact, computationally tractable outcomes.}


We will present a discussion on the computational aspects of our method in Section \ref{dis}. The paper concludes with Section \ref{con} where we summarise our findings, and present some ideas for future research.

\section{Call centres with two levels}\label{General}

In this section, we work with a downsized problem. Consider a call centre with two levels of callers and two levels of agents, where level 1 agents can take only level 1 calls and level 2 agents can serve both level 1 and 2 calls. We aim to compute the policy that minimises the total expected cost of abandonments for a given number of agents. This optimal policy is fully flexible in the sense that it takes the optimal state-dependent action whenever there is an arrival or service without assuming that there is any predefined functional form of the policy.

\subsection{Events requiring decision making}\label{DM}

In the call centre system, two types of events result in an opportunity to allocate a customer to an unoccupied server -- when there is a call arrival, and when there is a service completion \footnote{In a finite horizon system, where a policy is possibly non-stationary, it might conceivably be optimal to move a queued level 1 customer to a free level 2 server at a time which is not triggered by an arrival or service. We believe that the effect of neglecting such a possibility is likely to be negligible.}. When a call arrives, we need to decide which available agent to assign that call to, or we can choose to add the call to a queue. When an agent serves a call, we need to decide which call out of those waiting (if any) we should assign to this agent, or we can choose to reserve the agent for another call that might arrive later.

In order to reduce our action space, we make the assumption that if a level 1 call arrives or is in the queue and a level 1 agent is available, the optimal action is to assign the level 1 call to the available level 1 agent. This makes sense, since level 1 agents are dedicated to level 1 calls and there is nothing to be gained by making a level 1 customer wait when such an agent is available. 

When an arrival occurs, we need to make a decision only when
\begin{itemize}
\item the arrival is a level 1 call, all level 1 agents are busy and level 2 agents are available, or 
\item the arrival is a level 2 call and level 2 agents are available.
\end{itemize}
The possible actions in these cases are to assign the arriving call to the available level 2 agent or not. We are not ruling out the case that level 2 agents might be reserved for level 1 calls, and hence might not take an arriving level 2 call even when they are available, although we expect that this will not happen for most sensible parameter values.

Similarly, when there is a service, we need to make a decision only when 
\begin{itemize}
\item the service is completed by a level 2 agent and either only level 1 or only level 2 calls are waiting in the queue, or
\item the service is completed by a level 2 agent and both level 1 and level 2 calls are waiting in the queue.
\end{itemize}
The possible actions in the first case are to assign the waiting call to the level 2 agent who just completed the service or not, and the possible actions in the second case are to assign a level 1 or a level 2 call to the level 2 agent who just completed the service. A level 2 agent cannot reject both level 1 and level 2 calls. 

\subsection{The model}\label{gen_model}

We model the system as a continuous-time Markov chain. The states are denoted by $(a_1,a_2,a_q,b_2,b_q)$, where $a_1$, $a_2$ and $a_q$ are the number of level 1 calls being served by level 1 agents, being served by level 2 agents, and waiting in the queue, respectively, and $b_2$ and $b_q$ are the number of level 2 calls being served by level 2 agents, and waiting in the queue, respectively.

Depending on the state, seven possible events can occur - an arrival of each level of caller, an abandonment by each level of caller, and a service of each level of caller by an agent at the same level, or a service of a level 1 caller by a level 2 agent. Whenever these events fall into one of the categories described in Subsection \ref{DM}, we need to make a decision.

For $i \in \{1,2\}$, let $k_i$ be the total number of agents at level $i$. Therefore, when the state is $\textbf{s} = (a_1,a_2,a_q,b_2,b_q)$ there are $L_1(\mathbf{s}) = k_1 - a_1$ agents available at level 1 and $L_2(\mathbf{s}) = k_2 - a_2 - b_2$ agents available at level 2. Let $\ell$ be the total capacity of the system, implying that the total number of callers at any point in the system can be at most $\ell$, including callers from all levels waiting in the queue and being served.

The state space $S_g \coloneqq \{(a_1,a_2,a_q,b_2,b_q)\}$ is constrained such that
\begin{itemize}
    \item $a_1 \in \{0,1,\dots,k_1\}$,
    \item $a_2 \in \{0,1,\dots,k_2\}$,
    \item $b_2 \in \{0,1,\dots,k_2-a_2\}$,
    \item $a_1 + a_2 + a_q + b_2 + b_q \leq \ell$,
    \item if $a_1 < k_1$, then $a_q = 0$, and
    \item if $a_2 + b_2 < k_2$, then $a_q = 0$ or $b_q = 0$.
\end{itemize}
The second last constraint rules out the case where level 1 calls are in the queue and level 1 agents are available. The last constraint rules out the case where both level 1 and level 2 calls are in the queue while level 2 agents are available.

For all $\mathbf{s} \in S_g$, let $\mathcal{E}(\mathbf{s})$ denote the set of possible events that can occur when the system is in state $\mathbf{s}$. Let $e_{ai}$ denote the arrival of a level $i$ caller, $e_{si}$ denote the completion of a service by a level $i$ agent, and $e_{abi}$ denote an abandonment by a level $i$ caller. We use $\mathcal{A}(\mathbf{s},e)$ to denote the set of possible actions that can be taken when the system is in state $\mathbf{s}$ and an event $e$ occurs. Let $A_i$ denote the action that an incoming call of level $i$ is assigned to a level $i$ agent, and $A_q$ represent the action that an incoming call is added to the queue. 
Similarly, when there is a service completion by an agent at level $i$, let $S_{i}$ represent the action that a level $i$ call from the queue is assigned to that agent, and $S_{q}$ represent the action of leaving the agent unoccupied.

For many states and events, it may happen that the action set contains only one action and does not require any decision making. For example, if the system is in a state such that $L_1(\mathbf{s})>0$, then $\mathcal{A}(\mathbf{s},e_{a1}) $ contains only $A_1$, since an arriving level 1 call will always be assigned to an available level 1 agent.

For $i \in \{1,2\}$, let $\lambda_i$ and $\theta_i$ be the arrival and the abandonment rate of level $i$ calls, respectively. For $i \in \{1,2\}$, let $\mu_i$ be the service rate of level $i$ calls by level $i$ agents and let $\mu_1'$ be the service rate of level 1 calls by level 2 agents. Let $\gamma_{i}$ be the cost incurred when a level $i$ call abandons. For a specific policy, which defines the action that needs to be taken in every state, it is straightforward to construct a transition matrix $Q$ for the underlying continuous-time Markov chain in terms of the parameters $\lambda_i$, $\mu_i$, $\mu_1'$ and $\theta_i$. From this, it is also straightforward to derive the transition matrix for the embedded discrete-time jump chain. An example of this construction is provided in Appendix \ref{SubS:Transition}. 

We now incorporate these policy-dependent discrete-time transition matrices into a Markov decision process formulation.

Let 
\begin{equation}
    \kappa(\mathbf{s}) = (\lambda_1 + \lambda_2)\mathbbm{1}_{a_1+a_2+a_q+b_2+b_q < \ell} + \mu_1 a_1  + \mu_1' a_2 + \mu_2 b_2 + \theta_1 a_q + \theta_2 b_q
\end{equation}
so that $\kappa(\mathbf{s})$ is the total transition rate from state $\mathbf{s}.$
For $i \in \{1,2\}$ and $\textbf{s} \in S_g$, let 
\begin{equation}\label{eq:bellman}
    p_{ai} = \lambda_i\mathbbm{1}_{a_1+a_2+a_q+b_2+b_q < \ell}/\kappa(\mathbf{s})
\end{equation} 
be the probability that a level $i$ call arrives, 
\begin{equation}
    p_{s1} = \mu_1 a_1/\kappa(\mathbf{s})
\end{equation} be the probability that a level 1 call's service is completed by a level 1 agent, 
\begin{equation}
    p_{s1}' = \mu_1' a_2/\kappa(\mathbf{s})
\end{equation} 
be the probability that a level 1 call's service is completed by a level 2 agent, 
\begin{equation}
    p_{s2} = \mu_2 b_2/\kappa(\mathbf{s})
\end{equation} 
be the probability that a level 2 call's service is completed by a level 2 agent, 
\begin{equation}
    p_{ab1} = \theta_1 a_q/\kappa(\mathbf{s})
\end{equation} 
be the probability that a level 1 call abandons, and 
\begin{equation}\label{eq:bellman_end}
    p_{ab2} = \theta_2 b_q/\kappa(\mathbf{s})
\end{equation} 
be the probability that a level 2 call abandons.

\subsection{A discrete time approximation}\label{dta}

In this section we apply the dynamic programming method of Richard Bellman (\textcite{bellman1957dynamic}) to the discrete-time embedded chain defined by equations \eqref{eq:bellman} -- \eqref{eq:bellman_end} over a fixed total number $M$ of time steps. We use the parameter values given in Table \ref{GenEg1}, along with $\ell = 20$ and $M = 100$.
\begin{table}[ht]
    \centering
    \begin{tabular}{|c|c|c|}
    \hline
         $i$ & 1 & 2 \\ \hline
         $\lambda_i$ & 1 & 1 \\
         $\mu_i$ & 1/4 & 1/4 \\
         $\mu_i'$ & 1/4 & - \\
         $\theta_i$ & 1 & 2 \\
         $k_i$ & 5 & 5 \\
         $\gamma_{i}$ & 1 & 2 \\ \hline
    \end{tabular}
    \caption{Parameter values for the discrete-time example.}
    \label{GenEg1}
\end{table}  

For $m \in \{0,1,\dots,M\}$ and $\mathbf{s} \in S_g$, let $V(m,\textbf{s})$ be the minimum expected cost of abandonments between time steps $m$ and $M$ given that we are in state $\mathbf{s}$ at time step $m$. Since we only want to minimise the expected cost of abandonments \textit{between} time steps $m$ and $M$, we take $V(M,\textbf{s}) = 0$ $\forall \mathbf{s} \in S_g$.

Depending on the state $\textbf{s}$ at time step $m$, there will be certain actions that need to be taken whenever there is an event. 

Consider the state $\textbf{s}= (5,1,1,2,0)$. There are five level 1 callers being served by level 1 agents, one level 1 caller being served by a level 2 agent, one level 1 caller waiting in the queue, two level 2 callers being served by level 2 agents, two level 2 agents unoccupied and no level 2 callers waiting in the queue. In this state, there are two possible events that require a decision -- the arrival of a level 1 caller or the completion of a service by a level 2 agent. For each of these events $e$, there are two possible actions in the set $\mathcal{A}(\mathbf{s},e)$, resulting in a total of four possible sets of actions.

\begin{itemize}
    \item The first set of actions is to assign an arriving level 1 caller to an available level 2 agent, and to assign a level 1 caller waiting in the queue to a level 2 agent if there is a service completion by a level 2 agent. Under this set of actions,
    \begin{equation}
    \label{V1}
        \begin{split}
            V(m,(5,1,1,2,0)) &= p_{a1} V(m+1,(5,2,1,2,0)) + p_{a2} V(m+1,(5,1,1,3,0)) +\\& p_{s1} V(m+1,(5,1,0,2,0)) + p_{s1}' V(m+1,(5,1,0,2,0)) +\\& p_{s2} V(m+1,(5,2,0,1,0)) + p_{ab1} (V(m+1,(5,1,0,2,0))+\gamma_{1}),
        \end{split}
    \end{equation}
    
    \item The second set of actions is to send an arriving level 1 caller to the queue, and to assign a level 1 caller waiting in the queue to a level 2 agent if there is a service completion by a level 2 agent. Under this set of actions,
    \begin{equation}
    \label{V2}
        \begin{split}
             V(m,(5,1,1,2,0)) &= p_{a1} V(m+1,(5,1,2,2,0)) + p_{a2} V(m+1,(5,1,1,3,0)) +\\& p_{s1} V(m+1,(5,1,0,2,0)) + p_{s1}' V(m+1,(5,1,0,2,0)) +\\& p_{s2} V(m+1,(5,2,0,1,0)) + p_{ab1} (V(m+1,(5,1,0,2,0))+\gamma_{1}),
        \end{split}
    \end{equation}
    \item The third set of actions is to assign an arriving level 1 caller to an available level 2 agent, but not to assign a level 1 caller waiting in the queue to a level 2 agent when there is a service completion by a level 2 agent. Under this set of actions,
    \begin{equation}
    \label{V3}
        \begin{split}
             V(m,(5,1,1,2,0)) &= p_{a1} V(m+1,(5,2,1,2,0)) + p_{a2} V(m+1,(5,1,1,3,0)) +\\& p_{s1} V(m+1,(5,1,0,2,0)) + p_{s1}' V(m+1,(5,0,1,2,0)) +\\& p_{s2} V(m+1,(5,1,1,1,0)) + p_{ab1} (V(m+1,(5,1,0,2,0))+\gamma_{1}),
        \end{split}
    \end{equation}
    \item The fourth set of actions is to send an arriving level 1 caller to the queue, and not to assign a level 1 caller waiting in the queue to a level 2 agent when there is a service completion by a level 2 agent. Under this set of actions,
    \begin{equation}
    \label{V4}
        \begin{split}
             V(m,(5,1,1,2,0)) &= p_{a1} V(m+1,(5,1,2,2,0)) + p_{a2} V(m+1,(5,1,1,3,0)) +\\& p_{s1} V(m+1,(5,1,0,2,0)) + p_{s1}' V(m+1,(5,0,1,2,0)) +\\& p_{s2} V(m+1,(5,1,1,1,0)) + p_{ab1} (V(m+1,(5,1,0,2,0))+\gamma_{1}).
        \end{split}
    \end{equation}
    
\end{itemize}

For given values of $V(m+1,\textbf{s})$ on the right hand sides of equations \eqref{V1} -- \eqref{V4}, we select the set of actions that results in the minimum value of $V(m,(5,1,1,2,0))$. By proceeding backwards in this way from time step $M-1$ to $0$, and considering all states, we obtain the optimal action for each state and time step. Together, these optimal actions form the optimal policy. 
This policy can be described in an action matrix with the number of rows  equal to the size of the state space and the number of columns equal to $M+1$.

For the set of parameters given in Table \ref{GenEg1}, we first discuss the optimal policy at $m = 0$. We find that it is optimal to assign level 2 calls to level 2 agents whenever there is a level 2 agent available irrespective of the number of level 1 callers in the queue.

When all agents are busy and both level 1 and level 2 calls are waiting in the queue, it is optimal to assign a level 2 call to a level 2 agent who becomes free, and when only level 1 calls are in the queue, the optimal action is to leave all the level 1 calls in the queue level 1 irrespective of how many there are.

When all level 1 agents are busy and level 2 agents are available, there are two subcases -- when no one is waiting in the queue, we assign an incoming level 1 call to a level 2 agent only if there are more than two level 2 agents available, and when level 1 calls are waiting in the queue, the optimal policy is described in Table \ref{action1} ($L_2(\mathbf{s})$ is the number of level 2 agents available out of 5 and $a_q$ is the number of level 1 callers in the queue).  
\begin{table}[ht]
    \centering
    \begin{tabular}{|p{1cm}|p{1cm}|p{9cm}|}
    \hline
    $L_2(\mathbf{s})$ & $a_q$ & Optimal action\\\hline
    3-5 & 1 & Assign a level 1 call to a level 2 agent when there is an arrival or a service.\\
    2 & 1 & Do not assign a level 1 call to a level 2 agent when there is an arrival but assign a level 1 call to a level 2 agent when there is a service.\\
    1 & 1 & Do not assign a level 1 call to a level 2 agent in the case of an arrival or service.\\ \hline
    2-5 & 2 & Assign a level 1 call to a level 2 agent when there is an arrival or a service.\\
    1 & 2 & Do not assign a level 1 call to a level 2 agent in the case of an arrival or service.\\ \hline
    2-5 & 3-10 & Assign a level 1 call to a level 2 agent when there is an arrival or a service.\\
    1 & 3-10 & Do not assign level 1 call to a level 2 agent when there is an arrival but assign a level 1 call to a level 2 agent when there is a service.\\ \hline
    1-4 & 11-15 & Assign a level 1 call to a level 2 agent when there is an arrival or a service.\\ \hline
    \end{tabular}
    \caption{Optimal actions for the case when all the level 1 agents are busy, some level 2 agents are available and some level 1 calls are waiting in the queue.}
    \label{action1}
\end{table}

As can be observed, the policy becomes more and more flexible (or requires less reservation of level 2 servers) as the number of level 1 calls in the queue increases. When there is only one level 1 customer in the queue, we start with reserving two agents at level 2 for level 2 calls only, then we decrease the reservation to one and finally, we remove it completely. Note that in the last row of the table, although we have written $11-15$ for $a_q$, the number of callers in the queue can increase only to the level allowed by the maximum capacity. For instance, if all level 1 agents and three level 2 agents are busy, this implies that there are already $8$ callers in the system, so the number of callers in the queue can increase only to $12$.

As we get closer to the time horizon, the optimal policy becomes more flexible in allocating level 1 calls to level 2 servers. The policy remains as in Table \ref{action1} for the first $72$ steps, and then it changes slightly. After $87$ steps, it starts assigning level 1 calls to level 2 agents whenever there is more than one level 2 agent available, and for the last three time steps, it assigns level 1 calls to level 2 agents whenever there is any availability. 

As we can see in Table \ref{GenEg1}, the abandonment rate and the cost of abandonment for level 2 calls are both twice that of level 1 calls. If we consider examples where the difference in abandonment rates and costs between the two levels is not as significant, the optimal policy obtained from this method is to assign a level 1 call to a level 2 agent whenever a level 2 agent is available, and to prioritise assigning a level 2 call when both types of calls are present. Later in the paper, we will refer to such a policy as the {\it zero-reservation policy}.

The overall pattern of the optimal policy that we observe is that if the number of level 1 calls waiting in the queue is small, then we should reserve more level 2 agents exclusively for level 2 calls. As the number of level 1 calls waiting increases, we should assign level 1 calls to level 2 agents with more flexibility. The optimal policy also depends on how close the time horizon -- closer to the time horizon, the optimal policy becomes increasingly flexible in assigning level 1 callers to level 2 agents.

\subsection{An infinite time horizon Markov decision process with discounting}\label{cta}

In this section, in order to get some more insight into the type of policies that are optimal for our system, we find the optimal stationary policy for an approximate infinite horizon continuous-time Markov decision process with discounting  (\textcite{puterman2014markov}, \textcite{xie2019optimizing}, \textcite{baird1994reinforcement}).

Since we are now working with an infinite time horizon, the value function depends only on the state and not on the time. For a positive discounting factor $\delta$, let $C_\delta(\mathbf{s})$ be a random variable denoting the discounted number of abandonments over the infinite horizon given that we start in state $\mathbf{s} \in S_g$  and $V(\mathbf{s}) = \mathbbm{E}(C_\delta(\mathbf{s}))$. Let $T_0$ be the time to the next event, then $T_0$ is a random variable distributed exponentially with rate $\kappa(\mathbf{s})$. Now, we can write expressions for the expected cost of abandonments given the time to the next event $\mathbbm{E}(C_\mathbf{s}|T_0 = t)$ to calculate $V(\mathbf{s})$.

For example, again consider state $(5,1,1,2,0)$. We saw in Section \ref{dta} that, for this state, there are two possible actions when a level 1 call arrives, and two actions possible when a level 2 agent completes a service. Hence,
\begin{equation}
\begin{split}    
    \mathbbm{E}(C_\delta(5,1,1,2,0)|T_0 = t) =& e^{-\delta t}(p_{a1} \min(V(5,2,1,2,0),V(5,1,2,2,0)) +\\& p_{a2} V(5,1,1,3,0) + p_{s1} V(4,1,1,2,0) +\\& \min(p_{s1}' V(5,1,0,2,0) + p_{s2} V(5,2,0,1,0),\\&p_{s1}' V(5,0,1,2,0) + p_{s2} V(5,1,1,1,0)) +\\& p_{ab1} (V(5,1,0,2,0)+\gamma_1)).
\end{split}
\end{equation}

Now, using the independence of $T_0$ and the subsequent evolution of the process, we can write
\begin{equation}\label{eq:disc}
\begin{split}    
    V(5,1,1,2,0) =& \left(\frac{\kappa(5,1,1,2,0)}{\delta + \kappa(5,1,1,2,0)}\right)(p_{a1} \min(V(5,2,1,2,0),V(5,1,2,2,0)) \\& + p_{a2} V(5,1,1,3,0) + p_{s1} V(5,1,0,2,0) \\& +\min(p_{s1}' V(5,1,0,2,0) + p_{s2} V(5,2,0,1,0),\\&p_{s1}' V(5,0,1,2,0) + p_{s2} V(5,1,1,1,0)) \\& + p_{ab1} (V(5,1,0,2,0)+\gamma_1)).
\end{split}
\end{equation}
Similarly we can write equations for the value function in every state. These equations are non-linear because they contain $\min(\cdot)$. We use an iterative method to obtain an approximate solution (\textcite{powell2007approximate}, \textcite[Chapter 1]{chang2013simulation}). 

Let $K$ be a fixed number of iterations and $\alpha \in [0,1]$ be a step size. Then, for $k \in \{0,1,\dots,K-1\}$, we write 
\begin{equation}\label{eq:iter}
\begin{split}    
    V_{k+1}(5,1,1,2,0) =& \alpha V_{k}(5,1,1,2,0)\\& +(1-\alpha)\left[\left(\frac{\kappa(5,1,1,2,0)}{\delta +  \kappa(5,1,1,2,0)}\right)(p_{a1} \min(V_k(5,2,1,2,0),V_k(5,1,2,2,0))\right. \\& + p_{a2} V_k(5,1,1,3,0) + p_{s1} V_k(5,1,0,2,0)\\& +\min(p_{s1}' V_k(5,1,0,2,0) + p_{s2} V_k(5,2,0,1,0),\\&p_{s1}' V_k(5,0,1,2,0) + p_{s2} V_k(5,1,1,1,0)) \\& \left. + p_{ab1} (V_k(5,1,0,2,0)+\gamma_1))\right],
\end{split}
\end{equation}
with analogous equations for every other state.
We can initialise the iteration by choosing the $V_0(\mathbf{s})$ arbitrarily. The values of $\alpha$ and $V_0(\mathbf{s})$ will affect the convergence speed of the algorithm, but not the final values $V_K(\mathbf{s})$ and the optimal policy. We can also increase the computational efficiency of this method by introducing a stopping criterion that compares the vector of value functions at iteration $k+1$ with that at iteration $k$. We obtain the optimal stationary policy by extracting the arguments that lead to the minimum values on the right hand side.

\subsubsection{Example}

We tested the algorithm with two different initialisations of $V_0(\mathbf{s})$ -- one using the values from the model discussed in Section \ref{dta} and the other just simply putting them equal to zero. We also tested it with different values of $\alpha$. All combinations led to the same optimal policy and very close value functions with $K = 1000$. 

We used the above procedure to calculate the optimal policy, again using the parameter values in Table \ref{GenEg1}. Again, it is always optimal to assign level 2 calls to level 2 agents whenever there is an availability. When all the level 1 agents are busy while level 2 agents are available and no one is waiting in the queue, we assign an incoming level 1 call to a level 2 agent only if there are more than two level 2 agents available. In the case when level 1 calls are waiting in the queue, the policy is described in Table \ref{action2}.

\begin{table}[ht]
    \centering
    \begin{tabular}{|p{1cm}|p{1cm}|p{9cm}|}
    \hline
    $L_2(\mathbf{s})$ & $a_q$ & Optimal action\\\hline
    2-5 & 1 & Assign a level 1 call to a level 2 agent when there is an arrival or a service.\\
    1 & 1 & Do not assign a level 1 call to a level 2 agent in the case of an arrival or service.\\ \hline
    2-5 & 2-10 & Assign a level 1 call to a level 2 agent when there is an arrival or a service.\\
    1 & 2-10 & Do not assign level 1 call to a level 2 agent when there is an arrival but assign a level 1 call to a level 2 agent when there is a service.\\ \hline
    1-5 & 11-15 & Assign a level 1 call to a level 2 agent when there is an arrival or a service.\\ \hline
    \end{tabular}
    \caption{Optimal actions for the case when all the level 1 agents are busy, some level 2 agents are available, and some level 1 calls are waiting in the queue.}
    \label{action2}
\end{table}

When all agents are busy and both level 1 and level 2 calls are waiting in the queue, it is optimal to assign a level 2 call to a level 2 agent who becomes free. When only level 1 calls are in the queue, the optimal action is not to assign a level 1 call to a level 2 agent from the queue regardless of how many are waiting. We observe that the results are similar to those obtained from the discrete time model in Table \ref{action1}.

\section{Reservation policies}\label{ch:res}

We can characterise the optimal policies that we derived in Sections \ref{dta} and \ref{cta} as {\it bi-threshold policies} in which we take actions depending on the number of free agents available at level 2 and the number of level 1 calls waiting in the queue. These policies assign level 1 calls to level 2 servers if both the numbers of unoccupied level 2 agents and level 1 calls waiting in the queue are above respective thresholds. 

We also observe from the finite-horizon model in Section \ref{dta} that in finite time the optimal policy changes in favour of allocating level 1 calls to level 2 servers as the number of time steps left until the planning horizon decreases. 

It is reasonable to say that such a policy could be too complicated to implement in practice. It is computationally expensive as it requires a large number of state space variables, even for a two-level problem  and the exact policy is sensitive to changes in the parameter values. 

Based upon the intuition gained from looking at the models in Sections \ref{dta} and \ref{cta}, we propose a class of policies that are likely to be close to be optimal for hierarchical models with $\nu \geq 2$. We call these policies \textit{reservation policies}. They will form the focus of the rest of this paper.

\begin{definition}\label{Def}
A \textit{reservation policy} for a hierarchical model of the class defined in Section \ref{intro} is a policy under which calls are assigned to agents based on reservation thresholds $\Theta_i$ for each level $i$ agent. Under this policy

\begin{itemize}
\item if there is an arrival of a level $i-1$ caller,
\begin{itemize}
\item first, attempt to assign the caller to an available level $i-1$ agent,
\item if all level $i-1$ agents are busy, then check the number of available level $i$ agents,
\item if at least $\Theta_i$ level $i$ agents are available, assign the call to one of them,
\item if fewer than $\Theta_i$ level $i$ agents are available, assign the caller to the queue.
\end{itemize}
\item if a level $i$ agent completes a service,
\begin{itemize}
    \item the agent first checks for level $i$ callers in the queue and serves such a caller if any exist
    \item if there are no level $i$ callers in the queue, but there are level $i-1$ callers waiting -- the agent will serve a level $i-1$ caller if the number of available level $i$ agents (including the newly available agent) is at least $\Theta_i$.
\end{itemize}
\end{itemize}
\end{definition}

We see that a reservation policy depends only on the number of higher level agents available and not on the number of lower level callers callers in the queue. Such a policy is parameterised by the values $\mathbf{\Theta} = (\Theta_2,\ldots,\Theta_\nu)$ of the thresholds.

Using the method that we develop in Section \ref{ab4} for a model with $\nu=2$, we analytically calculated the expected cost of abandonments $C(\Theta_2,T)$ for the parameter values in Table \ref{GenEg1} and time duration $(0,t)$ for three different values of $t$ and all the possible values of $\Theta_2 \in \{0,1,\dots,k_2\}$. The values are given in Table \ref{lv2eg}.

\begin{table}[ht!]
\centering
\begin{subtable}{0.3\textwidth}
\centering
\begin{tabular}{|c|c|}
\hline
$\Theta_2$ & $C(\Theta_2,2)$ \\
\hline
0 & 0.0078 \\
1 & 0.0078 \\
2 & 0.0081 \\
3 & 0.0086 \\
4 & 0.0094 \\
5 & 0.0100 \\
\hline
\end{tabular}
\caption{$T = 2$}
\end{subtable}
\begin{subtable}{0.3\textwidth}
\centering
\begin{tabular}{|c|c|}
\hline
$\Theta_2$ & $C(\Theta_2,15)$ \\
\hline
0 & 4.467 \\
1 & 4.278 \\
2 & 4.290 \\
3 & 4.387 \\
4 & 4.476 \\
5 & 4.511 \\
\hline
\end{tabular}
\caption{$T = 15$}
\end{subtable}
\begin{subtable}{0.3\textwidth}
\centering
\begin{tabular}{|c|c|}
\hline
$\Theta_2$ & $C(\Theta_2,60)$ \\
\hline
0 & 27.071 \\
1 & 25.706 \\
2 & 25.528 \\
3 & 25.864 \\
4 & 26.204 \\
5 & 26.332 \\
\hline
\end{tabular}
\caption{$T = 60$}
\end{subtable}
\caption{Expected cost of abandonments for parameter values given in Table \ref{GenEg1}.}
    \label{lv2eg}
\end{table}

We observe that for small values of $T$ (analogous to being near the finite horizon in the discrete time method), the minimum expected cost of abandonments occurs at $\Theta_2 = 0$, that is, the optimal policy is the zero reservation policy where we do not reserve any level 2 agents and assign a level 1 call to a level 2 agent whenever there is an availability. For $T=15$ and values around that, the minimum is achieved at $\Theta_2 = 1$ implying that it is optimal to reserve one level 2 agent for level 2 calls only, and for the larger values of $T$, it is optimal to reserve two agents at level $2$. This is similar to what we observe in Section \ref{dta}. 

Overall, when we compare the optimal reservation policy with the optimal policies obtained in Sections \ref{dta} and \ref{cta}, we find that the differences arise primarily in scenarios where a large number of level 1 callers are waiting in the queue.

We implemented the policies derived in Section \ref{dta} at $m = 0$ and Section \ref{cta}, along with the reservation policy with $\Theta_2 = 2$ on a simulated underlying continuous-time Markov chain using the parameter values in Table \ref{GenEg1}. We observed that the expected cost of abandonments is nearly identical across all three policies and the performance was not distinguishable just on the basis of the simulation results.

Moreover, when the abandonment rates and the associated costs are similar across the two levels, the optimal reservation policy corresponds to a zero reservation policy, the same as the policy derived in Sections \ref{dta} and \ref{cta}.

In Section \ref{ab4}, we shall present a method to calculate the expected cost of abandonments and waiting time distributions of customers for a reservation policy in continuous time for a finite time horizon. Since reservation policies require a smaller number of variables compared to the policies discussed in Section \ref{General}, we are able to extend our analysis to the case where there are four levels of both agents and callers. 

\section{Performance measures for call centres with four levels} \label{ab4}

In this section, we will formulate a modification of the method introduced in \textcite{chiera_what_2002} for our industry partner's call centre which has four levels of calls and agents.
The objective will be to calculate the optimal thresholds $\mathbf{\Theta}$ and the corresponding expected costs of abandonments and waiting time distributions for the class of reservation policies defined in \ref{Def}.

We model the system as a continuous-time Markov chain, assuming that the calls are arriving according to Poisson processes, and the service and abandonment times are exponentially distributed.

For $i \in \{1,2,3,4\}$, let $\lambda_i$, $\theta_i$, and $\gamma_i$ be the arrival rate, abandonment rate, and abandonment cost of level $i$ callers, respectively. Let $k_i$ be the total number of level $i$ agents, $\mu_i$ be the service rate when level $i$ callers are served by level $i$ agents, and for $i \in \{1,2,3\}$, let $\mu_i'$ be the service rate when level $i$ callers are served by level $i+1$ agents. Let $\ell$ be the total capacity of the system shared by all callers regardless of their level and let $\beta$ be the cost incurred when customer is lost/blocked. Because the number of reserved agents has to be less than the total number of agents at that level, we have for $i \in \{2,3,4\}$, $\Theta_i \in \{0,1,\dots,k_i\}$.

\subsection{State space}\label{statespace}

We denote the state of the system by $(a,a_1,b,b_1,c,c_1,d)$, where $a, b, c$ and $d$ are the number of level 1, 2, 3 and 4 callers in the system, and $a_1, b_1$ and $c_1$ are the number of level 1, 2 and 3 callers in the system being served by level 2, 3 and 4 agents, respectively. This state space description is slightly different from that used in Section \ref{gen_model}, because the nature of the reservation policy allows us to convey the same amount of information using fewer variables.

When the system is in state $(a,a_1,b,b_1,c,c_1,d)$ there are $a$ level 1 callers, $a_1$ of which are being served by level 2 agents. Hence, $\min(a-a_1,k_1)$ level 1 callers are being served by level 1 agents, and $\max(0,a-(a_1+k_1))$ level 1 callers are waiting in the queue. It is possible that some of the level 1 agents are free while level 2 agents are serving level 1 callers if level 1 agents became free after the calls were already assigned to level 2 agents. Now, if $a_1$ level 2 agents are serving level 1 callers, and $b_1$ level 2 callers are being served by level 3 agents, then this implies $\min(b-b_1,k_2-a_1)$ level 2 callers are being served by level 2 agents, and $\max(0,b-(b_1+k_2-a_1))$ are waiting in the queue. Similarly, we can calculate how many level 3 and level 4 callers are waiting and how many are being served by each level of agent.

Here, $a,b,c,d \in \{0,1,\dots,\ell\}$ and $a+b+c+d \leq \ell$. The maximum value $a_1$ can take is $\min(k_2-\Theta_2,a)$ because by the definition of $a_1$, it cannot be more than the total number of level 1 callers and has to be less than or equal to the total number of level 2 agents. The minimum value of $a_1$ is zero in all cases except for when $a > k_1$ and $b < k_2-\Theta_2$ (that is, the number of level 1 callers are more than the number of level 1 agents and there are more than $\Theta_2$ level 2 agents available). Hence,\begin{equation}\min(a-k_1,k_2-b-\Theta_2) \leq a_1 \leq \min(k_2-\Theta_2,a).\end{equation} Similarly, \begin{equation}\min(b-k_2,k_3-c-\Theta_3)\mathbbm{1}_{b>k_2,c<k_3-\Theta_3} \leq b_1 \leq \min(k_3-\Theta_3,b)\end{equation} and \begin{equation}\min(c-k_3,k_4-d-\Theta_4)\mathbbm{1}_{c>k_3,d<k_4-\Theta_4} \leq c_1 \leq \min(k_4-\Theta_4,c).\end{equation}
Hence, the state space is given by $S = \{(a,a_1,b,b_1,c,c_1,d)\}$ such that
\begin{itemize}
    \item $a,b,c,d \in \{0,1,\dots,l\}$,
    \item $a+b+c+d \leq \ell$,
    \item $\min(a-k_1,k_2-b-\Theta_2) \mathbbm{1}_{a > k_1,b < k_2-\Theta_2} \leq a_1 \leq \min(k_2-\Theta_2,a)$, 
    \item $\min(b-k_2,k_3-c-\Theta_3)\mathbbm{1}_{b>k_2,c<k_3-\Theta_3} \leq b_1 \leq \min(k_3-\Theta_3,b)$ and 
    \item $\min(c-k_3,k_4-d-\Theta_4)\mathbbm{1}_{c>k_3,d<k_4-\Theta_4} \leq c_1 \leq \min(k_4-\Theta_4,c)$.
\end{itemize}

\subsection{Transition rates}
{Transitions between states occur when there is an arrival, a service or an abandonment at one of the levels 1 to 4. For example, given that the system is in state $(a,a_1,b,b_1,c,c_1,d)$, if a level 1 caller arrives, there are three possibilities.
    \begin{itemize}
        \item The caller is allocated to a level 1 agent if at least one of them is available leading to the state $(a+1,a_1,b,b_1,c,c_1,d)$.
        \item The caller is allocated to a level 2 agent (if no level 1 agent is available, no level 2 call is in waiting, and there are more than $\Theta_2$ level 2 agents available) leading to the state $(a+1,a_1+1,b,b_1,c,c_1,d)$.
        \item The caller joins the queue of level one customers if neither of the first two conditions are satisfied leading to the state $(a+1,a_1,b,b_1,c,c_1,d)$.
    \end{itemize} 
Hence, there is a transition to the state $(a+1,a_1+\mathbbm{1}(a-a_1 \geq k_1, b - b_1 < k_2 - \Theta_2 - a_1),b,b_1,c,c_1,d)$ with rate $\lambda_1$, if $a+b+c+d < \ell$, otherwise the caller is lost, costing $\beta$.
    
A complete list of all the transition rates for this model are given in Appendix \ref{SubS:Transition}.}

\subsection{The expected cost of losses and abandonments}\label{Subs:LT_Ab4}

For a given state $(a,a_1,b,b_1,c,c_1,d) \in S$ at time $0$ and reservation vector $\bm{\Theta}$, define $C_{a,a_1,b,b_1,c,c_1,d}(\bm{\Theta},T), \,\,$ $T \geq 0$ to be the expected cost of abandonments and losses during the time interval $(0,T)$.

Using a methodology similar to that used in \textcite{chiera_what_2002}, {we apply first-step analysis followed by the Laplace transform to compute $C_{a,a_1,b,b_1,c,c_1,d}(\mathbf{\Theta},t)$ for all possible values of $\mathbf{\Theta}$}. We condition on the first time that the system leaves its original state, writing an expression for $C_{a,a_1,b,b_1,c,c_1,d}(\mathbf{\Theta},T \mid x)$ where $x$ is the sojourn time in the state  $(a,a_1,b,b_1,c,c_1,d)$. 

If $T < x$, the system stays in its original state for the whole interval and there are no abandonments in time $(0,T)$. However, if the system is at its capacity, there might be lost customers. The expected cost in that case is given by $ T\, \beta\,\sum_{i = 1}^4 \lambda_i$. 

If $T\geq x$ and, for example, the system is currently in a state such that, $a < k_1, b+a_1 < k_2, c+b_1 < k_3, d+c_1 > k_4$, implying that there are agents available at each level except level 4, and there are customers waiting at level 4, then the expected cost of abandonments and losses is given by 

\begin{equation}\label{eq:4level_test}
    C_{a,a_1,b,b_1,c,c_1,d}(\bm{\Theta},T|x) = \frac{C^n_{a,a_1,b,b_1,c,c_1,d}(\bm{\Theta},T|x)}{\kappa_{a,a_1,b,b_1,c,c_1,d}(\bm{\Theta})}
\end{equation}
where
\begin{equation}
\label{eq:cond}
\begin{split}    
C^n_{a,a_1,b,b_1,c,c_1,d}(\bm{\Theta},T|x) = & \lambda_1 C_{a+1,a_1,b,b_1,c,c_1,d}(\bm{\Theta},T-x)\\& + \lambda_2 C_{a,a_1,b+1,b_1,c,c_1,d}(\bm{\Theta},T-x) \\& + \lambda_3 C_{a,a_1,b,b_1,c+1,c_1,d}(\bm{\Theta},T-x) \\& + \lambda_4 C_{a,a_1,b,b_1,c,c_1,d+1}(\bm{\Theta},T-x) \\&
+ (a-a_1)\mu_1 C_{a-1,a_1,b,b_1,c,c_1,d}(\bm{\Theta},T-x)\\& + (b-b_1)\mu_2 C_{a,a_1,b-1,b_1,c,c_1,d}(\bm{\Theta},T-x)\\& + (c-c_1)\mu_3 C_{a,a_1,b,b_1,c-1,c_1,d}(\bm{\Theta},T-x)\\& + (k_4 - c_1)\mu_4 C_{a,a_1,b,b_1,c,c_1,d-1}(\bm{\Theta},T-x)\\& 
+ a_1 \mu_1' C_{a-1,a_1-1,b,b_1,c,c_1,d}(\bm{\Theta},T-x)\\& + b_1 \mu_2' C_{a,a_1,b-1,b_1-1,c,c_1,d}(\bm{\Theta},T-x)\\& + c_1 \mu_3' C_{a,a_1,b,b_1,c-1,c_1-1,d}(\bm{\Theta},T-x)\\&
+ (d-(k_4-c_1))\theta_4(\gamma_4+C_{a,a_1,b,b_1,c,c_1,d-1}(\bm{\Theta},T-x))
\end{split}
\end{equation}
and
\begin{equation}
\begin{split}
    \kappa_{a,a_1,b,b_1,c,c_1,d}(\bm{\Theta}) = & \sum_{i = 1}^4\lambda_i + (a-a_1)\mu_1 + (b-b_1)\mu_2 + (c-c_1)\mu_3 + (k_4-c_1)\mu_4 \\& + a_1\mu_1' + b_1\mu_2' + c_1\mu_3' + (d-(k_4-c_1))\theta_4.
\end{split}
\end{equation}
If the system is at capacity, this equation needs to be modified to include the blocking cost in time $(0,x)$. After writing similar expressions for all states $C_{a,a_1,b,b_1,c,c_1,d}(\bm{\Theta},T|x)$, we can integrate to remove the conditioning on $x$, observing that $x$ is an exponentially distributed random variable. The right hand side of Equation \eqref{eq:cond} becomes a sum of convolution integrals which suggests that we should take Laplace transforms. This results in expressions of the form
\begin{equation}
\label{eq:LT}
\Tilde{C}_{a,a_1,b,b_1,c,c_1,d}(\bm{\Theta},s) = \dfrac{\xi^a_{a,a_1,b,b_1,c,c_1,d}(\bm{\Theta},s)}{\psi^a_{a,a_1,b,b_1,c,c_1,d}(s)}
\end{equation}
for the Laplace transform of the expected cost of losses and abandonments in time $(0,T)$, conditional on the system being in state $(a,a_1,b,b_1,c,c_1,d)$ at time 0. 

The exact expressions for $\xi_{a,a_1,b,b_1,c,c_1,d}(\mathbf{\Theta},s)$ and $\psi_{a,a_1,b,b_1,c,c_1,d}(s)$
are given in Appendix \ref{laplace}. 

We can solve equations \eqref{eq:LT} numerically and then invert them using the Euler method discussed in \textcite{abate1995numerical}.

If we let $\beta = 0$, $\gamma_i = 1$, and $\gamma_j = 0$ for $i \neq j$, we get an expression for the expected number of abandonments by level $i$ callers and if we let $\gamma_i = 0$ and $\beta = 1$ we get an expression for the expected number of lost customers.

{Using the same logic as above, we can also compute the time-dependent transient probabilities by formulating the Kolmogorov differential equations and applying the Laplace transform to them.}

We can also calculate the expected waiting time of customers, which is proportional to the expected cost of abandonments with $\beta = 0$. This is not surprising. While class $i$ customers are waiting, they each have their own `exponential abandonment process’ going, which has rate $\theta_i$ and is independent of what other customers are doing. We can also calculate the expected number of services and the expected waiting time of the $j^{\text{th}}$ customer in the queue.

\subsection{Numerical results}\label{num}

Before proceeding to the computation of waiting time distributions, we first present numerical examples to demonstrate how the proposed method can be applied to calculate the expected number of abandonments and, subsequently, to determine the optimal policy.

\subsubsection{Example 1}\label{eg1}

Consider a model with parameter values given in Table \ref{eg1t} along with $\beta = 0$, $\ell = 10$, and whose initial state is $\bm{\nu_0} = (0,0,0,0,0,0,0)$.

\begin{table}[ht]
    \centering
    \begin{tabular}{|c|c|c|c|c|}\hline
        $i$ & 1 & 2 & 3 & 4 \\ \hline
        $\lambda_i$ & 1 & 1/2 & 1/5 & 1/8 \\[0.2ex]
        $\mu_i$ & 2/3 & 1/2 & 1/4 & 1/6 \\ [0.2ex]
        $\mu_i'$ & 2/3 & 1/2 & 1/4 & - \\
        $\theta_i$ & 2 & 1 & 1 & 1 \\
        $k_i$ & 3 & 2 & 2 & 2\\
        $\gamma_i$ & 1 & 1 & 1 & 1\\
        \hline
    \end{tabular}
    \caption{Parameter values for Example 1.}
    \label{eg1t}
\end{table}

Assuming that the unit of time is one minute, the values in Table \ref{eg1t} imply that calls are arriving one per minute, one per two minutes, one per five minutes, and one per eight minutes for level $1$, $2$, $3$ and $4$, respectively, capturing the fact that more complex calls are rarer. The service rate for the four levels of calls is such that the expected service time increases with the increase in the level, and is indifferent to whether the call is being served by the same level agent or by a higher level agent. Similarly, the abandonment rate for level 1 calls is two per minute and for all other levels is one per minute. The cost of abandonment is equal for all four levels of calls and we are not counting the blocked customers (the expected number of callers blocked in this example is less than $0.3\%$ showing that $\ell = 10$ is reasonable). {The parameter values are chosen to reflect the characteristics of a real call centre.}

We calculated the expected number of call abandonments during $(0,60)$ minutes. Table \ref{t1} gives the expected cost of abandonments for all possible values of ${\bm{\Theta}}$.

\begin{table}[ht!]
\centering
\begin{subtable}{0.45\textwidth}
\centering
\begin{tabular}{|c|c|c|c|}
\hline
$\Theta_2$ & $\Theta_3$ & $\Theta_4$ & $C_{\bm{\nu_0}}(\bm{\Theta},60)$ \\
\hline
0 & 0 & 0 & 3.53 \\
0 & 0 & 1 & 4.03 \\
0 & 0 & 2 & 4.81 \\
0 & 1 & 0 & 5.22 \\
0 & 1 & 1 & 5.62 \\
0 & 1 & 2 & 6.21 \\
0 & 2 & 0 & 7.23 \\
0 & 2 & 1 & 7.57 \\
0 & 2 & 2 & 8.10 \\
1 & 0 & 0 & 5.37 \\
1 & 0 & 1 & 5.86 \\
1 & 0 & 2 & 6.59 \\
1 & 1 & 0 & 6.81 \\
1 & 1 & 1 & 7.19 \\
\hline
\end{tabular}
\end{subtable}
\hfill
\begin{subtable}{0.45\textwidth}
\centering
\begin{tabular}{|c|c|c|c|}
\hline
$\Theta_2$ & $\Theta_3$ & $\Theta_4$ & $C_{\bm{\nu_0}}(\bm{\Theta},60)$ \\
\hline
1 & 1 & 2 & 7.77 \\
1 & 2 & 0 & 8.51 \\
1 & 2 & 1 & 8.86 \\
1 & 2 & 2 & 9.40 \\
2 & 0 & 0 & 7.13 \\
2 & 0 & 1 & 7.61 \\
2 & 0 & 2 & 8.33 \\
2 & 1 & 0 & 8.45 \\
2 & 1 & 1 & 8.84 \\
2 & 1 & 2 & 9.42 \\
2 & 2 & 0 & 10.04 \\
2 & 2 & 1 & 10.38 \\
2 & 2 & 2 & 10.92 \\
\hline
\end{tabular}
\end{subtable}
\caption{Expected cost of abandonments for the parameter values given in Table \ref{eg1t}.}
\label{t1}
\end{table}

As can be seen from the table, the optimal value of ${\bm{\Theta}}$ is $(0,0,0)$, that is, not to use reservation of agents at all and allocate lower level calls to higher level agents as they are available. In this example, the optimal expected cost of abandonments is $3.53$. As a comparison, we can observe that, for the worst policy (when there is complete reservation), the expected cost of abandonments is $10.92$.

\subsubsection{Example 2}\label{eg2}

In this example, we look at a scenario in which level 4 calls are much more important than other calls, and if a level 4 caller abandons the queue it is four times more expensive than an abandonment by any other level of customer. {This setting could represent an airline call centre, where level 4 callers correspond to frequent business-class customers.}

The parameter values are given in Table \ref{eg2t} along with $\beta = 0$ and $\ell = 10$. Table \ref{t2} gives the expected cost of abandonments during $(0,60)$ minutes for all possible values of $\bm{\Theta}$.
\begin{table}[ht]
    \centering
    \begin{tabular}{|c|c|c|c|c|}\hline
        $i$ & 1 & 2 & 3 & 4 \\ \hline
        $\lambda_i$ & 1 & 1/2 & 1/4 & 1/4 \\[0.2ex]
        $\mu_i$ & 2/3 & 1/2 & 1/4 & 1/2 \\ [0.2ex]
        $\mu_i'$ & 2/3 & 1/2 & 1/4 & - \\
        $\theta_i$ & 2 & 1 & 1 & 1 \\
        $k_i$ & 3 & 2 & 2 & 2\\
        $\gamma_i$ & 1 & 1 & 1 & 4\\
        \hline
    \end{tabular}
    \caption{Parameter values for Example 2.}
    \label{eg2t}
\end{table}

\begin{table}[ht!]
\centering
\begin{subtable}{0.45\textwidth}
\centering
\begin{tabular}{|c|c|c|c|}
\hline
$\Theta_2$ & $\Theta_3$ & $\Theta_4$ & $C_{\bm{\nu_0}}(\bm{\Theta},60)$ \\
\hline
0 & 0 & 0 & 7.78 \\
0 & 0 & 1 & 7.13 \\
0 & 0 & 2 & 7.79 \\
0 & 1 & 0 & 9.13 \\
0 & 1 & 1 & 8.63 \\
0 & 1 & 2 & 9.17 \\
0 & 2 & 0 & 10.78 \\
0 & 2 & 1 & 10.28 \\
0 & 2 & 2 & 10.78 \\
1 & 0 & 0 & 9.53 \\
1 & 0 & 1 & 8.90 \\
1 & 0 & 2 & 9.53 \\
1 & 1 & 0 & 10.67 \\
1 & 1 & 1 & 10.15 \\
\hline
\end{tabular}
\end{subtable}
\hfill
\begin{subtable}{0.45\textwidth}
\centering
\begin{tabular}{|c|c|c|c|}
\hline
$\Theta_2$ & $\Theta_3$ & $\Theta_4$ & $C_{\bm{\nu_0}}(\bm{\Theta},60)$ \\
\hline
1 & 1 & 2 & 10.69 \\
1 & 2 & 0 & 12.07 \\
1 & 2 & 1 & 11.58 \\
1 & 2 & 2 & 12.06 \\
2 & 0 & 0 & 11.26 \\
2 & 0 & 1 & 10.63 \\
2 & 0 & 2 & 11.26 \\
2 & 1 & 0 & 12.30 \\
2 & 1 & 1 & 11.80 \\
2 & 1 & 2 & 12.32 \\
2 & 2 & 0 & 13.60 \\
2 & 2 & 1 & 13.10 \\
2 & 2 & 2 & 13.60 \\
\hline
\end{tabular}
\end{subtable}
\caption{Expected cost of abandonments for the parameter values given in Table \ref{eg2t}.}
\label{t2}
\end{table}

The expected cost of abandonments, in this case, ranges from $7.13$ to $13.60$, and the best reservation policy is $\bm{\Theta} = (0,0,1)$, that is, out of the two level 4 agents, we should reserve one for just level 4 callers and all other agents should be flexible.

\subsection{Waiting time distributions}

{We now extend the method described in Section~\ref{Subs:LT_Ab4} to compute the waiting time distributions of customers. The objective is to determine the probability that a call arriving within the interval $(0, T)$ is answered within a specified time frame, for instance, 20 seconds, given the current state of the system. This enables performance-based staffing and routing decisions. For instance, if the client requires that at least $80\%$ of calls be answered within $20$ seconds, we can determine the staffing levels and call allocation policy such that the probability of answering a call within $20$ seconds is at least $0.8$.}

\subsubsection{Expected time spent in each state during the time interval \texorpdfstring{$\bm{(0,T)}$}{(0,T)}}

{We start by calculating the expected time that the Markov chain spends in each state during the interval $(0,T)$, as this will later be used to calculate the waiting time distributions of customers.

We model the system as a continuous-time Markov chain and use the same state space as defined in Subsection~\ref{statespace} and transition rates as described in Subsection~\ref{SubS:Transition}.

For $T>0$ and reservation vector $\bm{\Theta}$, let $\bm{v}_{a,a_1,b,b_1,c,c_1,d}(\bm{\Theta},T)$ be a vector whose entries denote the expected time spent by the system in each state during the time interval $(0,T)$, given that the system is initially in state \((a,a_1,b,b_1,c,c_1,d) \in S\) at time \(0\).

As in Subsection~\ref{Subs:LT_Ab4}, we condition on the first time that the system leaves its original state, writing an expression for $\bm{v}_{a,a_1,b,b_1,c,c_1,d}(\bm{\Theta},T \mid x)$, where $x$ is the sojourn time in state $(a,a_1,b,b_1,c,c_1,d)$.

Let $\bm{e}_{a,a_1,b,b_1,c,c_1,d}$ be a unit vector with a $1$ in the entry corresponding to state $({a,a_1,b,b_1,c,c_1,d})$ for $({a,a_1,b,b_1,c,c_1,d}) \in S$.

Now if $T<x$, the system spends the entire time in the initial state, hence, $\bm{v}_{a,a_1,b,b_1,c,c_1,d}(\bm{\Theta},T \mid x) = T \bm{e}_{a,a_1,b,b_1,c,c_1,d}$.

If $T>x$, then the system spends $x$ time units in state $({a,a_1,b,b_1,c,c_1,d})$, and then transitions to a different state according to the transition matrix. The equation is derived in the same way as Equation \eqref{eq:4level_test}.

We can construct equations for $\bm{v}_{a,a_1,b,b_1,c,c_1,d}(\bm{\Theta},T)$ by integrating $\bm{v}_{a,a_1,b,b_1,c,c_1,d}(\bm{\Theta},T \mid x)$ with respect to $x$. We then take the Laplace transform of these equations. The set of linear equations satisfied by the Laplace transform of $\bm{v}_{a,a_1,b,b_1,c,c_1,d}(\bm{\Theta},T)$, denoted by ${\Tilde{\bm{v}}}_{a,a_1,b,b_1,c,c_1,d}(\bm{\Theta},s)$, is provided in Appendix~\ref{Ts1}.}

\subsubsection[Distribution of the state of the system observed by a uniformly distributed arrival in \texorpdfstring{${(0,T)}$}{(0,T)}]{Distribution of the state of the system observed by a uniformly distributed arrival in \texorpdfstring{$\bm{(0,T)}$}{(0,T)}}\label{dist4}

{The component corresponding to state $\bm{s} \in S$ of the vector $\bm{v}_{a,a_1,b,b_1,c,c_1,d}(\bm{\Theta},T)$ represents the expected amount of time that the system spends in state $\bm{s}$ during the interval $(0,T)$ . Consequently, $\bm{v}_{a,a_1,b,b_1,c,c_1,d}(\bm{\Theta},T)/T$ is a vector whose component corresponding to state $\bm{s}$ gives the expected proportion of time in $(0,T)$ during which the system is in state $\bm{s}$.

Now, assume that a tagged caller arrives at a uniformly distributed time $U$ within $(0,T)$. Then, the component corresponding to state $\bm{s}$ of $\bm{v}_{a,a_1,b,b_1,c,c_1,d}(\bm{\Theta},T)/T$ equals the expected value of the indicator random variable that the caller finds the system in state $\bm{s}$. This is, in fact, the probability that the caller finds the system in state $\bm{s}$. We denote this probability by $r_{(a,a_1,b,b_1,c,c_1,d),\bm{s}}(T)$.

We will use these probabilities to calculate the waiting time distributions. However, they can also be used to compute other important performance measures for call centres, such as the average queue size or the queue size distribution at a uniformly distributed point in $(0,T)$. For example, we might be interested in determining the minimum number of agents required such that the probability that the number of customers waiting in the system, observed at a uniformly distributed arrival time, exceeds a given threshold value $q$ is less than $0.05$. Let $B$ be an event when the number of customers waiting in the system exceeds $q$. Using the notation introduced above, this probability can be expressed as
\[
\sum_{i \in B} r_{(a,a_1,b,b_1,c,c_1,d),i}(T),
\]
and the condition becomes
\[
\sum_{i \in B} r_{(a,a_1,b,b_1,c,c_1,d),i}(T) < 0.05.
\]

Next, we calculate the probability that a caller is served within $y$ time units, given the state of the system at their arrival, that is, the conditional waiting time distribution of a tagged caller. This probability differs for each caller level. We begin with level~4 callers, as this is the simplest case, and then proceed to discuss the other levels in turn.}

\subsubsection{The waiting time distribution for a level 4 caller}

We can model the system for a tagged level~$4$ caller as a continuous-time Markov chain with state space 
\[
\{(c_1,d):\, c_1 + d \leq \ell,\; c_1 \leq k_4-\Theta_4\},
\]
where $d$ denotes the number of level~4 callers in front of the tagged caller. For a level~$4$ caller, the number of level~$1$ or $2$ callers in the system is irrelevant. Its waiting time distribution depends only on the number of level~$3$ callers being served by level~$4$ agents at the time of arrival. As soon as a level~$4$ agent becomes free, the arriving level~$4$ caller will have priority over callers at all other levels.

For $y > 0$, let $\Xi^4_{c_1,d}(\bm{\Theta},y)$ denote the probability that a level~4 caller is served within $y$ time units, given that the state of the system at their arrival is $(c_1,d)$ and the reservation vector is $\bm{\Theta}$. Furthermore, let $\tilde{\Xi}_{c_1,d}^4(\bm{\Theta},s)$ be its Laplace transform.

The transition rates of concern are those corresponding to the service of a level~3 caller by a level~4 agent, the service of a level~4 caller, and the abandonment of all level~4 callers including and ahead of the tagged customer in the queue.

The Laplace transform $\tilde{\Xi}_{c_1,d}^4(\bm{\Theta},s)$ equals $0$ for $c_1+d = \ell$ (since if the system is full, the arriving caller is blocked and therefore the probability of being served is $0$), and for $c_1+d < \ell$ we have
{\begin{equation}\small
\tilde{\Xi}_{c_1,d}^4(\bm{\Theta},s) = \begin{cases}
    \dfrac{1}{s} & d < k_4 - c_1 \\[2ex]
    \dfrac{c_1\mu_3'\tilde{\Xi}_{c_1-1,d}^4(\bm{\Theta},s) + \left((k_4-c_1)\mu_4 + (d-(k_4-c_1))\theta_4 \right)\tilde{\Xi}_{c_1,d-1}^4(\bm{\Theta},s)}{c_1\mu_3' + (k_4-c_1)\mu_4 + (d+1-(k_4-c_1))\theta_4 + s} & d \geq k_4 - c_1.
\end{cases}\end{equation}}

\subsubsection{The waiting time distribution for a level 3 caller}\label{lev3}

We can model the system for a tagged level~$3$ caller as a continuous-time Markov chain with state space
{\footnotesize\begin{equation*}
\big\{(b_1,c,c_1,d):\; b_1 + c + d \leq \ell,\; b_1 \leq k_3-\Theta_3,\; 
\min(c-k_3,\,k_4-d-\Theta_4)\mathbbm{1}_{\{c>k_3,\,d<k_4-\Theta_4\}} \leq c_1 \leq \min(k_4-\Theta_4,c)\big\},
\end{equation*}}

\noindent where $c$ denotes the number of level~3 callers in front of the tagged caller and $d$ denotes all level~4 callers. The waiting time of a level~$3$ caller depends on the number of level~$3$ callers ahead, the total number of level~$4$ callers in the system, and the number of level~$2$ callers being served by level~$3$ agents. 

Furthermore, while the arrival of new level~$3$ callers does not change the tagged caller's waiting time, the arrival of new level~$4$ callers is relevant because level~$4$ callers have priority over level~$3$ callers when being served by level~$4$ agents. The transition rates will be defined accordingly.

For $y > 0$, let ${\Xi}_{b_1,c,c_1,d}^3(\bm{\Theta},y)$ denote the probability that a level~$3$ caller is served within $y$ time units, given that the state of the system at their arrival is $(b_1,c,c_1,d)$ and the reservation vector is $\bm{\Theta}$.

The Laplace transforms, $\tilde{\Xi}_{b_1,c,c_1,d}^3(\bm{\Theta},s) = 0$ for $b_1 + c + d = \ell$ and otherwise
\begin{equation}\tilde{\Xi}_{b_1,c,c_1,d}^3(\bm{\Theta},s) = \begin{cases}
    \dfrac{1}{s}, & k_3 - b_1 > c - c_1 \text{ or } k_4 > d + c_1 + \Theta_4 \\[3ex]
    \dfrac{\xi^3_{b_1,c,c_1,d}(\bm{\Theta},s)}{\psi^3_{b_1,c,c_1,d}(s)}, & k_3 - b_1 \leq c - c_1  \text{ and } k_4 \leq d + c_1 + \Theta_4
\end{cases}\end{equation}
where
\begin{equation}
    \begin{split}
        \xi^3_{b_1,c,c_1,d}(\bm{\Theta},s) = & \lambda_4 \mathbbm{1}_{b_1+c+d < \ell-1} \tilde{\Xi}_{b_1,c,c_1,d+1}^3 + (k_3-b_1)\mu_3 \tilde{\Xi}_{b_1,c-1,c_1,d}^3(\bm{\Theta},s) + \\ & b_1\mu_2' \tilde{\Xi}_{b_1-1,c,c_1,d}^3(\bm{\Theta},s) + \\ & c_1\mu_3'\tilde{\Xi}_{b_1,c-1,c_1-(1-\mathbbm{1}_{c+1-c_1 > k_3-b_1, d = k_4 - c_1 - \Theta_4}),d}^3(\bm{\Theta},s) + \\ & \min(k_4-c_1,d)\mu_4 \tilde{\Xi}_{b_1,c,c_1 + \mathbbm{1}_{c+1-c_1 > k_3-b_1, d = k_4-c_1-\Theta_4},d-1}^3(\bm{\Theta},s) + \\ & (c-(k_3-b_1+c_1))\theta_3 \tilde{\Xi}_{b_1,c-1,c_1,d}^3(\bm{\Theta},s) + \\ & \max(0,d-(k_4-c_1))\theta_4 \tilde{\Xi}_{b_1,c,c_1,d-1}^3(\bm{\Theta},s) 
    \end{split}
\end{equation}
and
\begin{equation}
    \begin{split}
        \psi^3_{b_1,c,c_1,d}(s) = &\lambda_4 \mathbbm{1}_{b_1+c+d < \ell-1} + (k_3-b_1)\mu_3 + b_1\mu_2' + c_1\mu_3' + \min(d,k_4-c_1)\mu_4 + \\ & (c+1-(k_3-b_1+c_1))\theta_3 + \max(0,d-(k_4-c_1))\theta_4 + s.
    \end{split}
\end{equation}

\subsubsection{The waiting time distribution for a level 2 caller}
The expression for the probability of being served within $y$ time units for a level~$2$ caller, denoted by ${\Xi}_{a_1,b,b_1,c,c_1,d}^2(\bm{\Theta},y)$, is similar to that in Subsection~\ref{lev3}, except that two additional variables are included in the state space. This is because we now consider all level~$2$, $3$, and $4$ callers, in addition to level~$1$ callers being served by level~$2$ agents. Accordingly, the relevant transition rates are included as well. 

The state space is given by $\{(a_1,b,b_1,c,c_1,d)\}$
such that
\begin{itemize}
    \item $a_1 + b + c + d \leq \ell$, 
    \item $a_1 \leq k_2-\Theta_2$,
    \item $\min(b - k_2,\, k_3 - c - \Theta_3)\mathbbm{1}_{\{b > k_2,\, c < k_3 - \Theta_3\}} \leq b_1 \leq \min(b,\, k_3 - \Theta_3)$, and
    \item $\min(c - k_3,\, k_4 - d - \Theta_4)\mathbbm{1}_{\{c > k_3,\, d < k_4 - \Theta_4\}} \leq c_1 \leq \min(c,\, k_4 - \Theta_4)$,
\end{itemize}
where $b$ denotes the number of level~$2$ callers ahead of the tagged caller, and $c$ and $d$ denote the total number of level~$3$ and level~$4$ callers in the system, respectively. 

The expressions for the Laplace transforms are provided in Appendix~\ref{l2ap}.

\subsubsection{The waiting time distribution for a level 1 caller}

For a tagged level~$1$ caller, the state space is the same as $S$. The transition rates are also similar to those in Subsection~\ref{SubS:Transition}}, except for some minor modifications, for example, the transition rate associated with the arrival of new level~$1$ callers is not included. Let ${\Xi}_{a,a_1,b,b_1,c,c_1,d}^1(\bm{\Theta},y)$ denote the probability that a level~$1$ caller is served within $y$ time units, given that the state of the system at their arrival is $(a,a_1,b,b_1,c,c_1,d)$ and the reservation vector is $\bm{\Theta}$. The expressions for the Laplace transforms of these probabilities are provided in Appendix~\ref{l1ap}.

We can solve the Laplace transform equations for all four levels of callers and then invert them using the Euler method to obtain the required probabilities.

\subsubsection{The waiting time distribution of any caller}

{We have now established how to compute the waiting time distribution of a tagged caller, conditional on the system state at their arrival, as well as the state distribution for a customer arriving uniformly over the interval $(0, T)$.

For each level of caller, we can apply the law of total probability along with the Poisson Arrivals See Time Averages (PASTA) property to calculate the unconditional probability that a caller arriving uniformly in $(0,T)$ is served within $y$ time units by conditioning on the state of the system at the time of arrival.}

For a set of reservation parameters $\bm{\Theta}$ and $j \in \{1,2,3,4\}$, let $P_{\bm{\nu_0}}^j(\bm{\Theta},y,T)$ denote the probability that a level~$j$ caller who arrives during the interval $(0,T)$ is served within $y$ time units, given that the initial state of the system at time $0$ is $\bm{\nu_0} = (a_0,a_{10},b_0,b_{01},c_0,c_{01},d_0)$. Hence, for level~1 callers,
\begin{equation}P_{\bm{\nu_0}}^1(\bm{\Theta},y,T) = \sum_{(a,a_1,b,b_1,c,c_1,d) \in S} \Xi_{a,a_1,b,b_1,c,c_1,d}^1(\bm{\Theta},y)r_{\bm{\nu_0},(a,a_1,b,b_1,c,c_1,d)}(\bm{\Theta},T),\end{equation}
for level 2 callers
\begin{equation}P_{\bm{\nu_0}}^2(\bm{\Theta},y,T) = \sum_{(a,a_1,b,b_1,c,c_1,d) \in S} \Xi_{a_1,b,b_1,c,c_1,d}^2(\bm{\Theta},y)r_{\bm{\nu_0},(a,a_1,b,b_1,c,c_1,d)}(\bm{\Theta},T),\end{equation}
for level 3 callers
\begin{equation}P_{\bm{\nu_0}}^3(\bm{\Theta},y,T) = \sum_{(a,a_1,b,b_1,c,c_1,d) \in S} \Xi_{b_1,c,c_1,d}^3(\bm{\Theta},y)r_{\bm{\nu_0},(a,a_1,b,b_1,c,c_1,d)}(\bm{\Theta},T),\end{equation}
and for level 4 callers
\begin{equation}P_{\bm{\nu_0}}^4(\bm{\Theta},y,T) = \sum_{(a,a_1,b,b_1,c,c_1,d) \in S} \Xi_{c_1,d}^4(\bm{\Theta},y)r_{{\bm{\nu_0}},(a,a_1,b,b_1,c,c_1,d)}(\bm{\Theta},T).\end{equation}

Now, by further conditioning on the type of arriving caller, we obtain the overall probability that a caller at any level is served within $y$ time units. Let $P_{\bm{\nu_0}}(\bm{\Theta},y,T)$ denote the probability that a caller at any level who arrives during the interval $(0,T)$ is served within $y$ time units, given that the state of the system at time $0$ is $\bm{\nu_0} = (a_0,a_{10},b_0,b_{01},c_0,c_{01},d_0)$. Hence,
\begin{equation}
P_{\bm{\nu_0}}(\bm{\Theta},y,T) = \frac{\displaystyle\sum_{j=1}^4 \lambda_j P_{\bm{\nu_0}}^j(\bm{\Theta},y,T)}{\sum_{j=1}^4 \lambda_j}
\end{equation}
since $\dfrac{\lambda_i}{\sum_{j=1}^4 \lambda_j}$ is the probability that an arriving customer is of level $i$.

\subsection{Numerical results}\label{num_2}

We now demonstrate how this method can be applied to compute the distribution function and to find the optimal policy for a given set of parameters.

\subsubsection{Distribution function}

Consider the parameter values in Table~\ref{tab_df}, along with $\ell = 10$ and the initial state $\bm{\nu_0} = (0,0,0,0,0,0,0)$.

\begin{table}[ht]
    \centering
    \begin{tabular}{|c|c|c|c|c|}\hline
        $j$ & 1 & 2 & 3 & 4 \\ \hline
        $\lambda_j$ & 1 & 1/2 & 1/4 & 1/6 \\[0.2ex]
        $\mu_j$ & 2/3 & 1/2 & 1/3 & 1/6 \\ [0.2ex]
        $\mu_j'$ & 2/3 & 1/2 & 1/3 & - \\
        $\theta_j$ & 1 & 1/2 & 1/4 & 1/4 \\
        $k_j$ & 2 & 2 & 1 & 1\\
        \hline
    \end{tabular}
    \caption{Parameter values for the distribution function example.}
    \label{tab_df}
\end{table}

\begin{figure}[ht!]
    \centering
    \begin{subfigure}{0.45\linewidth}
        \centering
        \includegraphics[width=\linewidth]{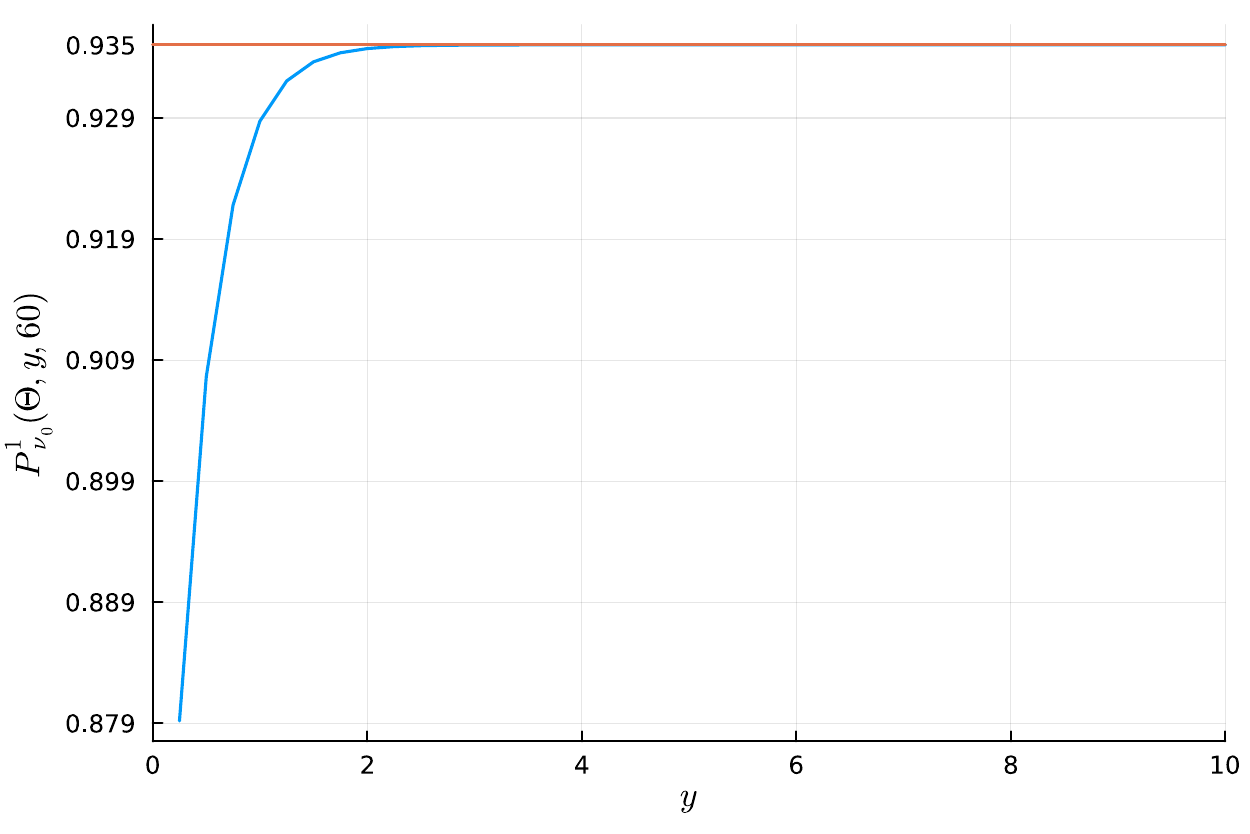}
        \caption{$P^1_{\nu_0}(\bm{\Theta},y,60)$ versus $y$}
        \label{fig_df1}
    \end{subfigure}
    \hfill
    \begin{subfigure}{0.45\linewidth}
        \centering
        \includegraphics[width=\linewidth]{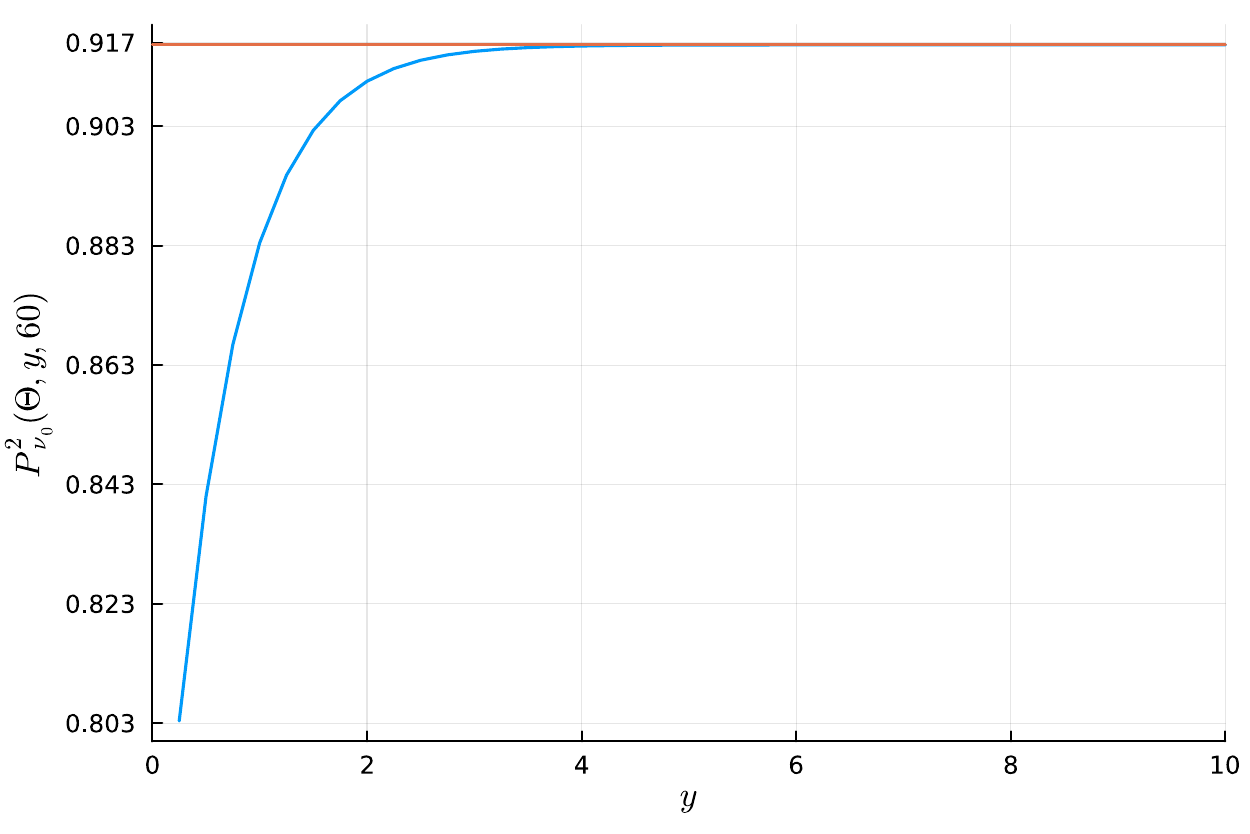}
        \caption{$P^2_{\nu_0}(\bm{\Theta},y,60)$ versus $y$}
        \label{fig_df2}
    \end{subfigure}
    
    \vskip\baselineskip
    \begin{subfigure}{0.45\linewidth}
        \centering
        \includegraphics[width=\linewidth]{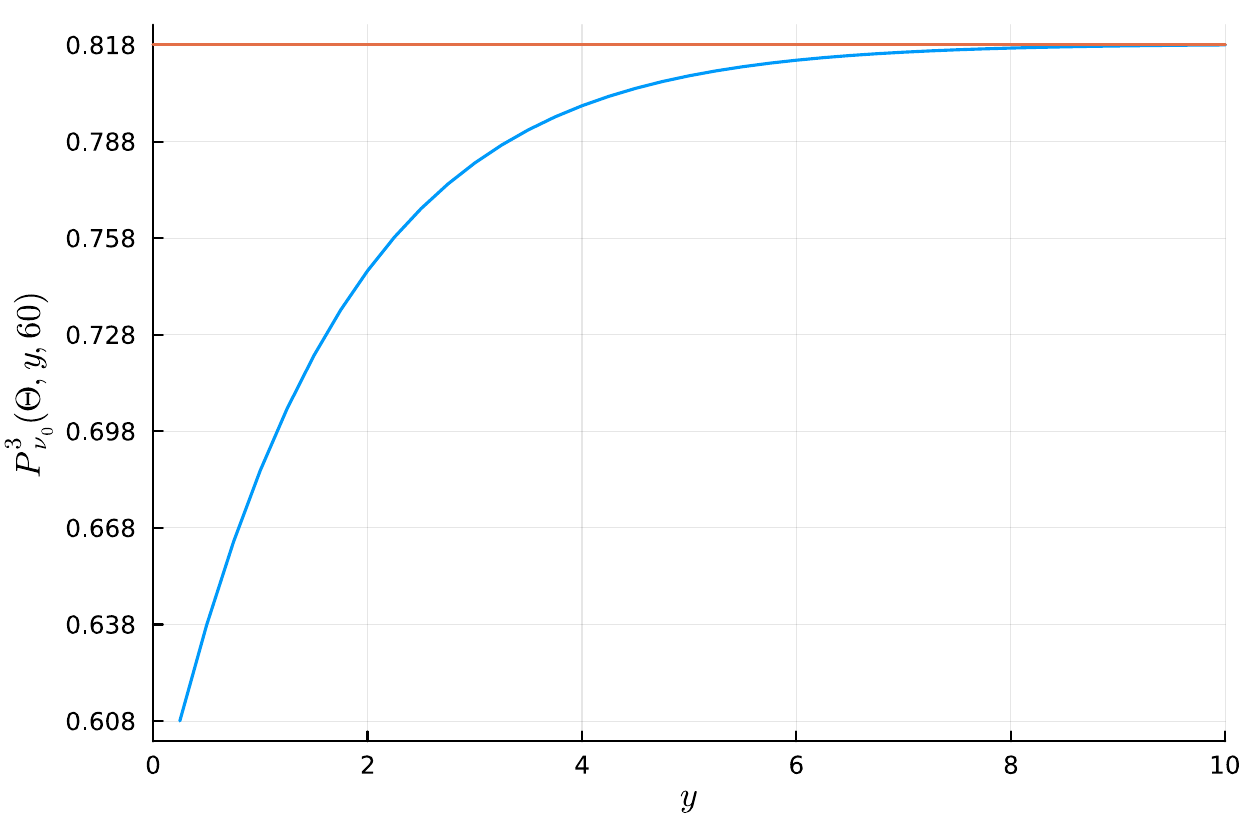}
        \caption{$P^3_{\nu_0}(\bm{\Theta},y,60)$ versus $y$}
        \label{fig_df3}
    \end{subfigure}
    \hfill
    \begin{subfigure}{0.45\linewidth}
        \centering
        \includegraphics[width=\linewidth]{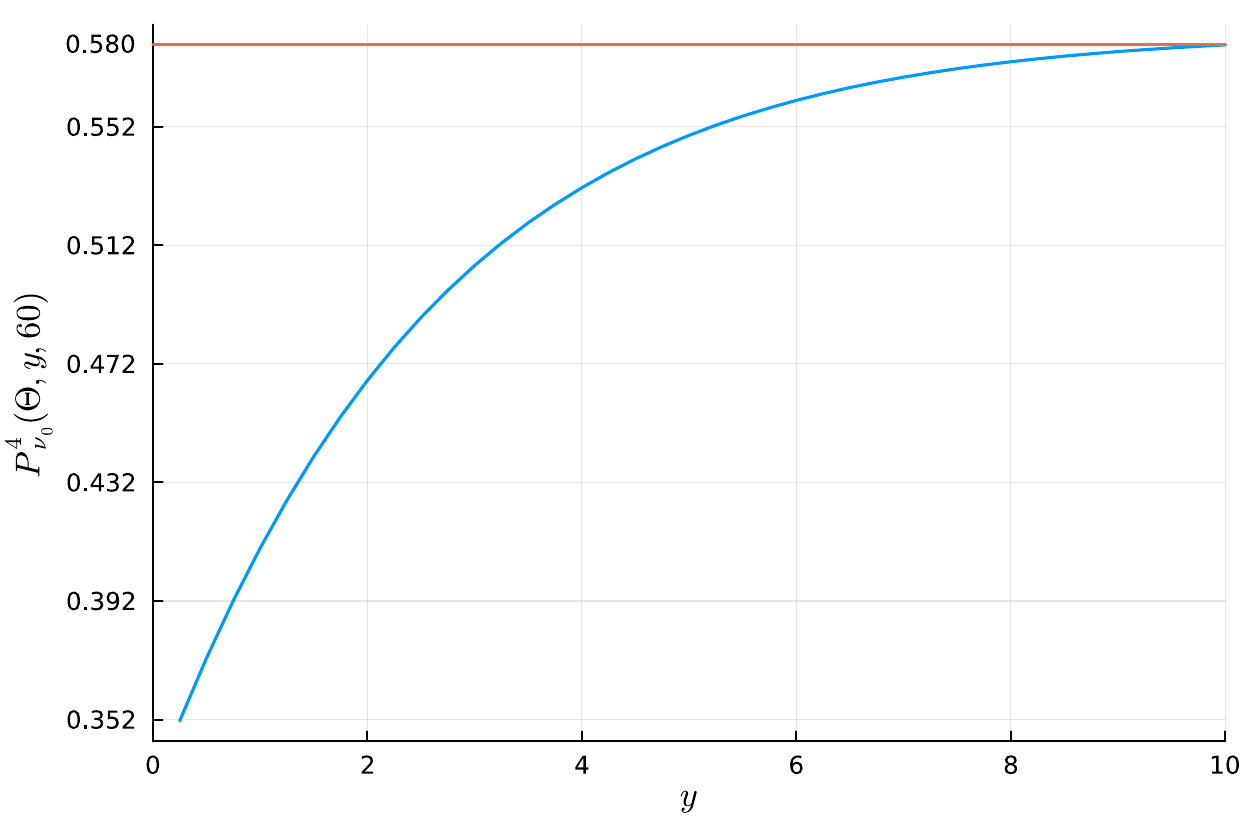}
        \caption{$P^4_{\nu_0}(\bm{\Theta},y,60)$ versus $y$}
        \label{fig_df4}
    \end{subfigure}

    \caption{Graphs of $P^k_{\nu_0}(\bm{\Theta},y,60)$ versus $y$ for $k \in {1,2,3,4}$, with parameter values in Table \ref{tab_df}, $\bm{\nu_0} = (0,0,0,0,0,0,0)$ and $\bm{\Theta} = (0,0,0)$.}
    \label{fig_df_all}
\end{figure}

\begin{figure}[ht!]
    \centering
    \includegraphics[width=0.5\linewidth]{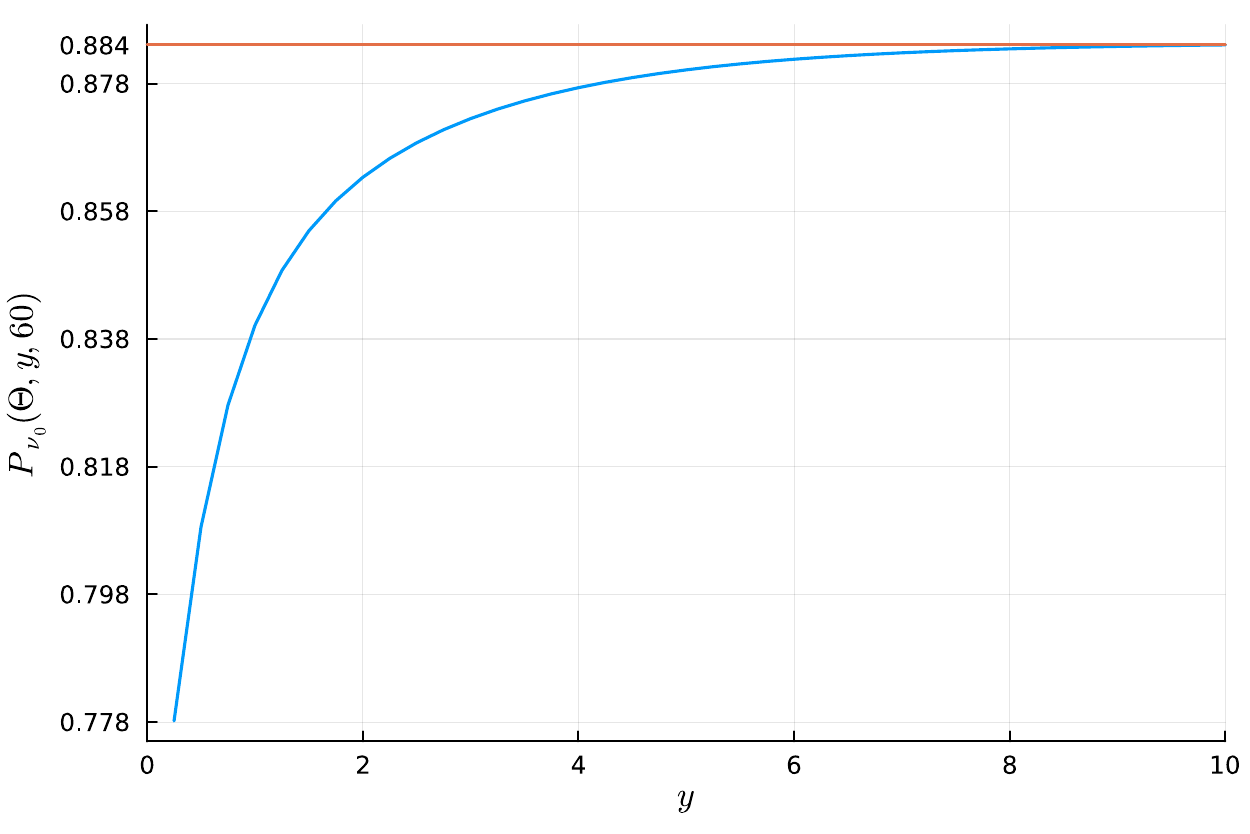}
    \caption{Graph of $P_{\nu_0}(\bm{\Theta},y,60)$ versus $y$ for parameter values in Table \ref{tab_df}, $\bm{\nu_0} = (0,0,0,0,0,0,0)$ and $\bm{\Theta} = (0,0,0)$.}
    \label{fig_df}
\end{figure}

Figures \ref{fig_df_all} and \ref{fig_df} present the distribution functions corresponding to the parameter values in Table~\ref{tab_df}, with Figure \ref{fig_df_all} illustrating the case for callers at each level and Figure \ref{fig_df} showing the general case under the zero-reservation policy.

In Figure \ref{fig_df_all}, we observe clear differences in the probability of being served within a given time period across the four levels. With $\bm{\Theta} = (0,0,0)$, level 1 callers have the highest probability of being served, while level 4 callers have the lowest. Furthermore, all curves eventually converge to the probability of ever being served, which is the complement of the probability of abandoning. For example, from Figure \ref{fig_df_all}, the abandonment probability for level 1 callers is $0.065$, whereas for level 4 callers it is $0.420$.

In Figure \ref{fig_df}, we can observe that the probability of serving a customer arriving uniformly in the next $60$ time units within $1$ time unit for the general case is approximately $0.83$. The probability of eventually serving a customer arriving in the next $60$ time units is approximately $0.88$, implying that the probability of a customer being lost or abandoning is approximately $0.12$.

\subsubsection{Example 1}

Now, we consider another example with parameter values given in Table~\ref{weg2t} along with $\ell = 10$. Let the initial state be $\bm{\nu_0} = (0,0,0,0,0,0,0)$. The parameter values are chosen such that customers are served faster when assigned to an agent at the same level than when assigned to an agent at a higher level.

\begin{table}[ht]
    \centering
    \begin{tabular}{|c|c|c|c|c|}\hline
        $j$ & 1 & 2 & 3 & 4 \\ \hline
        $\lambda_j$ & 1 & 1/2 & 1/4 & 1/8 \\[0.2ex]
        $\mu_j$ & 2/3 & 1/2 & 1/4 & 1/4 \\ [0.2ex]
        $\mu_j'$ & 2/3 & 1/8 & 1/16 & - \\
        $\theta_j$ & 2 & 1 & 1 & 1 \\
        $k_j$ & 3 & 2 & 1 & 1\\
        \hline
    \end{tabular}
    \caption{Parameter values for Example 2.}
    \label{weg2t}
\end{table}

We again calculate $P_{\bm{\nu_0}}(\bm{\Theta}, 1/3, 60)$. Table~\ref{t2w} presents the probabilities for the general case, while Tables~\ref{t21} and \ref{t22} show the probabilities at each level for all possible values of $\bm{\Theta}$.

As can be seen from Table~\ref{t2w}, the optimal value of $\bm{\Theta}$ for the general case is $(0,1,1)$, meaning that one agent is reserved each at levels $3$ and $4$. This is expected because, with the parameters given in Table~\ref{weg2t}, agents at levels $3$ and $4$ are less efficient in handling lower level queries. In this example, for the best reservation vector, the probability that a customer at any level is served within $20$ seconds is $0.832$, whereas in the worst-case reservation setting with $\bm{\Theta} = (2,0,0)$, this probability is $0.795$. Similar observations are found in Tables~\ref{t21} and \ref{t22}.

\begin{table}[ht!]
\centering
\begin{tabular}{|c|c|c|c|}
\hline
$\Theta_2$ & $\Theta_3$ & $\Theta_4$ & $P_{\bm{\nu_0}}(\bm{\Theta},1/3,60)$ \\
\hline
0 & 0 & 0 & 0.820 \\
0 & 0 & 1 & 0.829 \\
0 & 1 & 0 & 0.830 \\
0 & 1 & 1 & 0.832 \\
1 & 0 & 0 & 0.809 \\
1 & 0 & 1 & 0.818 \\
1 & 1 & 0 & 0.817 \\
1 & 1 & 1 & 0.819 \\
2 & 0 & 0 & 0.795 \\
2 & 0 & 1 & 0.803 \\
2 & 1 & 0 & 0.801 \\
2 & 1 & 1 & 0.803 \\
\hline
\end{tabular}
\caption{Probability that a customer arriving in the next $60$ minutes will be served within $20$ seconds for the parameter values given in Table \ref{weg2t}.}
\label{t2w}
\end{table}

\begin{table}[ht!]
\centering
\begin{subtable}{0.45\textwidth}
\centering
\begin{tabular}{|c|c|c|c|}
\hline
$\Theta_2$ & $\Theta_3$ & $\Theta_4$ & $P_{\bm{\nu_0}}^1(\bm{\Theta}, 1/3, 60)$ \\
\hline
0 & 0 & 0 & 0.959 \\
0 & 0 & 1 & 0.959 \\
0 & 1 & 0 & 0.964 \\
0 & 1 & 1 & 0.964 \\
1 & 0 & 0 & 0.922 \\
1 & 0 & 1 & 0.922 \\
1 & 1 & 0 & 0.925 \\
1 & 1 & 1 & 0.925 \\
2 & 0 & 0 & 0.889 \\
2 & 0 & 1 & 0.890 \\
2 & 1 & 0 & 0.889 \\
2 & 1 & 1 & 0.890 \\
\hline
\end{tabular}
\end{subtable}
\hfill
\begin{subtable}{0.45\textwidth}
\centering
\begin{tabular}{|c|c|c|c|}
\hline
$\Theta_2$ & $\Theta_3$ & $\Theta_4$ & $P_{\bm{\nu_0}}^2(\bm{\Theta}, 1/3, 60)$ \\
\hline
0 & 0 & 0 & 0.819 \\
0 & 0 & 1 & 0.838 \\
0 & 1 & 0 & 0.772 \\
0 & 1 & 1 & 0.772 \\
1 & 0 & 0 & 0.845 \\
1 & 0 & 1 & 0.861 \\
1 & 1 & 0 & 0.801 \\
1 & 1 & 1 & 0.801 \\
2 & 0 & 0 & 0.853 \\
2 & 0 & 1 & 0.869 \\
2 & 1 & 0 & 0.812 \\
2 & 1 & 1 & 0.812 \\
\hline
\end{tabular}
\end{subtable}
\caption{Probability that a level 1 customer (left) and a level 2 customer (right) arriving in the next $60$ minutes will be served within $20$ seconds for the parameter values given in Table \ref{weg2t}.}
\label{t21}
\end{table}

\begin{table}[ht!]
\centering
\begin{subtable}{0.45\textwidth}
\centering
\begin{tabular}{|c|c|c|c|}
\hline
$\Theta_2$ & $\Theta_3$ & $\Theta_4$ & $P_{\bm{\nu_0}}^3(\bm{\Theta}, 1/3, 60)$ \\
\hline
0 & 0 & 0 & 0.521 \\
0 & 0 & 1 & 0.366 \\
0 & 1 & 0 & 0.647 \\
0 & 1 & 1 & 0.501 \\
1 & 0 & 0 & 0.536 \\
1 & 0 & 1 & 0.381 \\
1 & 1 & 0 & 0.647 \\
1 & 1 & 1 & 0.501 \\
2 & 0 & 0 & 0.542 \\
2 & 0 & 1 & 0.387 \\
2 & 1 & 0 & 0.647 \\
2 & 1 & 1 & 0.501 \\
\hline
\end{tabular}
\end{subtable}
\hfill
\begin{subtable}{0.45\textwidth}
\centering
\begin{tabular}{|c|c|c|c|}
\hline
$\Theta_2$ & $\Theta_3$ & $\Theta_4$ & $P_{\bm{\nu_0}}^4(\bm{\Theta}, 1/3, 60)$ \\
\hline
0 & 0 & 0 & 0.309 \\
0 & 0 & 1 & 0.682 \\
0 & 1 & 0 & 0.359 \\
0 & 1 & 1 & 0.682 \\
1 & 0 & 0 & 0.314 \\
1 & 0 & 1 & 0.682 \\
1 & 1 & 0 & 0.359 \\
1 & 1 & 1 & 0.682 \\
2 & 0 & 0 & 0.316 \\
2 & 0 & 1 & 0.682 \\
2 & 1 & 0 & 0.359 \\
2 & 1 & 1 & 0.682 \\
\hline
\end{tabular}
\end{subtable}
\caption{Probability that a level 3 customer (left) and a level 4 customer (right) arriving in the next $60$ minutes will be served within $20$ seconds for the parameter values given in Table \ref{weg2t}.}
\label{t22}
\end{table}

Similarly, the optimal value of the reservation vector can be determined for any set of parameters, any given initial state, and any $y, T > 0$. 

Consider the parameter values in Table \ref{weg2t}, but this time with initial state  $\bm{\nu_0} = (1,0,2,0,1,0,0)$. We again calculate $P_{\bm{\nu_0}}(\bm{\Theta}, 1/3, 60)$. Table~\ref{t3w} presents the probabilities for the general case. As expected, the probability for all the threshold values has reduced in comparison to Table \ref{t2w} because there are already some customers at the beginning and the optimal threshold policy is still the same.

\begin{table}[ht!]
\centering
\begin{tabular}{|c|c|c|c|}
\hline
$\Theta_2$ & $\Theta_3$ & $\Theta_4$ & $P_{\bm{\nu_0}}(\bm{\Theta},1/3,60)$ \\
\hline
0 & 0 & 0 & 0.808 \\
0 & 0 & 1 & 0.818 \\
0 & 1 & 0 & 0.819 \\
0 & 1 & 1 & 0.821 \\
1 & 0 & 0 & 0.797 \\
1 & 0 & 1 & 0.806 \\
1 & 1 & 0 & 0.806 \\
1 & 1 & 1 & 0.808 \\
2 & 0 & 0 & 0.783 \\
2 & 0 & 1 & 0.792 \\
2 & 1 & 0 & 0.790 \\
2 & 1 & 1 & 0.792 \\
\hline
\end{tabular}
\caption{Probability that a customer arriving in the next $60$ minutes will be served within $20$ seconds for the parameter values given in Table \ref{weg2t} with $\bm{\nu_0} = (1,0,2,0,1,0,0)$.}
\label{t3w}
\end{table}

\section{Discussion}\label{dis}

{The methods presented so far have been discussed in the context of a single time interval. However, in practical applications, these methods would be implemented sequentially over multiple intervals. The analysis is equally valid for any initial state. In a real-world call centre, one could start with an empty system, determine the optimal policy for a given time interval, and then reapply the method at the end of that interval using the system's current state as the new initial state. This approach captures the time dependency of the optimal policy.

Even with the simplifications introduced by considering reservation policies instead of general policies, our analysis is constrained by the size of the system. As the system size increases, the matrix inversion required to calculate the Laplace transforms becomes larger, leading to higher computational costs. However, unlike when using Bellman's equation (Section \ref{dta}), changing $T$ does not increase the computational expense. The time to find the optimal policy also depends on the number of possible values of $\bm{\Theta}$.

For the type of reservation policy considered, the case with four levels of calls required seven state space variables. This was computationally tractable for relatively small values of the total system capacity $\ell$. For $\ell = 10$, computing the expected cost of abandonments for a given $\bm{\Theta}$ took approximately $300$ seconds for the hierarchical model with four levels (examples in Section \ref{num}). 

With efficient state space ordering, we have demonstrated feasibility up to $\ell = 15$. Memory allocation issues arose when $\ell > 15$ was used for the same model. For cases with fewer than four levels of calls (or policies with fewer than seven state space variables), this method can handle much larger call centres.

We made use of parallelisation in our computations wherever possible. We used The University of Melbourne's high-performing computer Spartan to run our code simultaneously for different reservation parameters using the array environment. In the Euler algorithm for the numerical inversion of Laplace transforms, there is a step where the transforms must be evaluated $30$ times at different values of $s$. For this step, we employed the \textit{Distributed} package in \textit{Julia}, which utilises multiple processors to divide the task and compile the results. When solving the inverse of the matrices for Laplace transforms, the use of sparse matrices and faster inversion methods has improved computational speed. These techniques have a greater impact when calculating distribution functions compared to expected abandonment costs, as the former are significantly more computationally expensive. For the examples presented in this Section \ref{num_2}, the final tables are produced in under $25$ minutes.

We believe that the method can be extended to five levels of calls. However, since each additional level requires adding two variables, and the state space grows exponentially with each variable, it would be challenging to extend the method beyond that point.



\section{Conclusion}\label{con}

{In this paper, we have employed a range of methods to identify the optimal finite-horizon call allocation policy for small call centres. We began by using backward induction and an infinite-horizon model with discounting to determine the class of policies that are likely to be optimal.

We then applied first-step analysis, Laplace transforms, and numerical inversions to derive two performance measures that are highly relevant to call centre operations -- the expected number of abandonments and the waiting time distributions of customers.

Dynamic programming and related methods allow us to perform policy optimisation, but they are not well-suited to the finite-horizon continuous-time setting. Conversely, the Laplace transform method allows precise calculation of performance measures in this setting for a given policy, but it does not support policy optimisation. In this paper, we combined these two approaches so that they complement one another, resulting in a complete framework for finding and analysing optimal call allocation in call centres.

Our analysis showed that reservation policies provide a relatively robust, adaptable, and computationally efficient alternative to fully flexible optimal policies, achieving near-optimal performance in call centres. Moreover, the developed methods are sufficiently general to be applied to other call centre systems and policies, provided the state space remains computationally tractable.}


\textbf{Acknowledgements.} This research was funded by the Australian Government through the Australian Research Council Industrial Transformation Training Centre in Optimisation Technologies, Integrated Methodologies, and Applications (OPTIMA), Project ID IC200100009. It was also supported by The University of Melbourne’s Research Computing Services and the Petascale Campus Initiative. The second author is also funded by the Melbourne Research Scholarship.

\textbf{Competing interests declaration.} The authors declare none.

\printbibliography[title= {References}]

\appendix

\section{Transition rates for the hierarchical model with  4 levels}\label{SubS:Transition}

Transitions between states occur when there is an arrival, a service or an abandonment at one of the levels 1 to 4. Below is a list of the transitions that ensue when each of these events occurs when the system is in state $(a,a_1,b,b_1,c,c_1,d)$. 

\begin{itemize}
    \item If a level 1 caller arrives, there are three possibilities.
    \begin{itemize}
        \item The caller is allocated to a level 1 agent if at least one of them is available leading to the state $(a+1,a_1,b,b_1,c,c_1,d)$.
        \item The caller is allocated to a level 2 agent (if no level 1 agent is available, no level 2 call is in waiting, and there are more than $\Theta_2$ level 2 agents available) leading to the state $(a+1,a_1+1,b,b_1,c,c_1,d)$.
        \item The caller joins the queue of level one customers if neither of the first two conditions are satisfied leading to the state $(a+1,a_1,b,b_1,c,c_1,d)$.
    \end{itemize} 
    Hence, there is a transition to the state $(a+1,a_1+\mathbbm{1}(a-a_1 \geq k_1, b - b_1 < k_2 - \Theta_2 - a_1),b,b_1,c,c_1,d)$ with rate $\lambda_1$ if $a+b+c+d < \ell$, otherwise the caller is lost, costing $\beta$.
    \item Similarly, a transition to the state $(a,a_1,b+1,b_1+\mathbbm{1}(b-b_1 \geq k_2, c - c_1 < k_3 - \Theta_3 - b_1),c,c_1,d)$ happens with rate $\lambda_2$ and a transition to the state $(a,a_1,b,b_1,c+1,c_1+\mathbbm{1}(c-c_1 \geq k_3, d < k_4 - \Theta_4 - c_1),d)$ happens with the rate $\lambda_3$ if $a+b+c+d < \ell$.
    \item When a level 4 caller arrives, there is only one possibility, that is, a transition to the state $(a,a_1,b,b_1,c,c_1,d+1)$ with rate $\lambda_4$ if $a+b+c+d < \ell$.
    \item A service of level 1 caller by a level 1 agent happens with rate $\min(k_1,a-a_1)\mu_1$ leading to the state $(a-1,a_1,b,b_1,c,c_1,d)$.
    \item When there is a service of a level 1 caller by a level 2 agent, there are two possibilities. 
    \begin{itemize}
        \item That agent becomes free or starts serving a level 2 caller if there are any in waiting leading to the state $(a-1,a_1-1,b,b_1,c,c_1,d)$.
        \item The agent starts serving another level 1 caller leading to the state $(a-1,a_1,b,b_1,c,c_1,d)$.
    \end{itemize}
    Hence, there is transition to the state $(a-1,a_1 - (1-\mathbbm{1}(a-a_1 > k_1, b - b_1 = k_2 - \Theta_2 - a_1),b,b_1,c,c_1,d)$ with rate $a_1 \mu_1'$.
    \item Similarly, when there is a service of a level 2 caller by a level 2 agent, that agent can either serve another level 2 caller, serve a level 1 caller or stay available. The system transitions to the state $(a,a_1+\mathbbm{1}(a-a_1 > k_1, b - b_1 = k_2 - \Theta_2 - a_1),b-1,b_1,c,c_1,d)$ with rate $\min(b-b_1,k_2-a_1)\mu_2$.
    \item Similarly, the system transitions to 
    \begin{itemize}
        \item $(a,a_1,b-1,b_1 - (1-\mathbbm{1}(b-b_1 > k_2 - a_1, c = k_3 - \Theta_3 - b_1),c,c_1,d)$ with rate $b_1 \mu_2'$ when there is a service of level 2 caller by a level 3 agent,
        \item $(a,a_1,b,b_1+\mathbbm{1}(b-b_1 > k_2, c - c_1 = k_3 - \Theta_3 - b_1),c-1,c_1,d)$ with rate $\min(c-c_1,k_3-b_1)\mu_3$ when there is a service of a level 3 caller by a level 3 agent,
        \item $(a,a_1,b,b_1,c-1,c_1 - (1-\mathbbm{1}(c-c_1 > k_3 - b_1, d = k_4 - \Theta_4 - c_1),d)$ with rate $c_1 \mu_3'$ when a level 3 caller is served by a level 4 agent and
        \item $(a,a_1,b,b_1,c,c_1+\mathbbm{1}(c-c_1 > k_3 - b_1, d = k_4 - \Theta_4 - c_1),d-1)$ with rate $\min(d,k_4-c_1)\mu_4$ when there is a service of level 4 caller.
    \end{itemize}
    \item The transitions in case of abandonments are straightforward. There is a transition to 
    \begin{itemize}
        \item $(a-1,a_1,b,b_1,c,c_1,d)$ when a level 1 caller abandons with rate $\max(0,a-(k_1 + a_1))\theta_1$ costing $\gamma_1$,
        \item $(a,a_1,b-1,b_1,c,c_1,d)$ when a level 2 caller abandons with rate $\max(0,b-(k_2-a_1+b_1))\theta_2$ costing $\gamma_2$,
        \item $(a,a_1,b,b_1,c-1,c_1,d)$ when a level 3 caller abandons with rate $\max(0,c-(k_3-b_1+c_1))\theta_3$ costing $\gamma_3$ and 
        \item $(a,a_1,b,b_1,c,c_1,d-1)$ when a level 4 caller abandons with rate $\max(0,d-(k_4-c_1))\theta_4$ costing $\gamma_4$.
    \end{itemize}
\end{itemize}

\newpage
\section{Expected cost of losses and abandonments for the hierarchical model with four levels}\label{laplace}

The Laplace transform
\begin{equation}
\Tilde{C}_{a,a_1,b,b_1,c,c_1,d}(\mathbf{\Theta},s) = \dfrac{\xi_{a,a_1,b,b_1,c,c_1,d}(\mathbf{\Theta},s)}{\psi_{a,a_1,b,b_1,c,c_1,d}(s)}
\end{equation}
where
\begin{equation}
    \begin{split}
        \xi_{a,a_1,b,b_1,c,c_1,d}(\mathbf{\Theta},s) = &(\lambda_1 \Tilde{C}_{a+1,a_1+\mathbbm{1}(a-a_1 \geq k_1, b - b_1 < k_2 - \Theta_2 - a_1),b,b_1,c,c_1,d}(\mathbf{\Theta},s) \\&+ \lambda_2 \Tilde{C}_{a,a_1,b+1,b_1+\mathbbm{1}(b-b_1 \geq k_2, c < k_3 - \Theta_3 - b_1),c,c_1,d}(\mathbf{\Theta},s)\\&+ \lambda_3 \Tilde{C}_{a,a_1,b,b_1,c+1,c_1+\mathbbm{1}(c-c_1 \geq k_3, d < k_4 - \Theta_4 - c_1),d}(\mathbf{\Theta},s)\\ & + \lambda_4 \Tilde{C}_{a,a_1,b,b_1,c,c_1,d+1}(\mathbf{\Theta},s)) \mathbbm{1}_{a+b+c+d < \ell}+\min(k_1,a-a_1)\mu_1 \Tilde{C}_{a-1,a_1,b,b_1,c,c_1,d}(\mathbf{\Theta},s) \\ &+ a_1 \mu_1' \Tilde{C}_{a-1,a_1 - (1-\mathbbm{1}(a-a_1 > k_1, b - b_1 = k_2 - \Theta_2 - a_1),b,b_1,c,c_1,d}(\mathbf{\Theta},s))\\ &+\min(b-b_1,k_2-a_1)\mu_2 \Tilde{C}_{a,a_1+\mathbbm{1}(a-a_1 > k_1, b - b_1 = k_2 - \Theta_2 - a_1),b-1,b_1,c,c_1,d}(\mathbf{\Theta},s)\\ & + b_1 \mu_2' \Tilde{C}_{a,a_1,b-1,b_1 - (1-\mathbbm{1}(b-b_1 > k_2 - a_1, c = k_3 - \Theta_3 - b_1),c,c_1,d}(\mathbf{\Theta},s)) \\&+\min(c-c_1,k_3-b_1)\mu_3 \Tilde{C}_{a,a_1,b,b_1+\mathbbm{1}(b-b_1 > k_2, c - c_1 = k_3 - \Theta_3 - b_1),c-1,c_1,d}(\mathbf{\Theta},s)\\ & + c_1 \mu_3' \Tilde{C}_{a,a_1,b,b_1,c-1,c_1 - (1-\mathbbm{1}(c-c_1 > k_3 - b_1, d = k_4 - \Theta_4 - c_1),d}(\mathbf{\Theta},s)) \\ &+ \min(d,k_4-c_1)\mu_4 \Tilde{C}_{a,a_1,b,b_1,c,c_1+\mathbbm{1}(c-c_1 > k_3 - b_1, d = k_4 - \Theta_4 - c_1),d-1}(\mathbf{\Theta},s) \\ &+\max(0,a-(k_1 + a_1))\theta_1(\Tilde{C}_{a-1,a_1,b,b_1,c,c_1,d}(\mathbf{\Theta},s) + \gamma_1/s)\\ &+\max(0,b-(k_2-a_1+b_1))\theta_2(\Tilde{C}_{a,a_1,b-1,b_1,c,c_1,d}(\mathbf{\Theta},s)+\gamma_2/s)\\ &+\max(0,c-(k_3-b_1+c_1))\theta_3(\Tilde{C}_{a,a_1,b,b_1,c-1,c_1,d}(\mathbf{\Theta},s)+\gamma_3/s)\\ &+\max(0,d-(k_4-c_1))\theta_4(\Tilde{C}_{a,a_1,b,b_1,c,c_1,d-1}(\mathbf{\Theta},s)+\gamma_4/s)\\ &+(\lambda_1+\lambda_2+\lambda_3+\lambda_4)\beta\mathbbm{1}_{a+b+c+d = \ell}/s,
    \end{split}
\end{equation}
\begin{equation}
    \begin{split}
        \psi_{a,a_1,b,b_1,c,c_1,d}(s) = &(\lambda_1+\lambda_2+\lambda_3+\lambda_4)\mathbbm{1}_{a+b+c+d < \ell} + \min(k_1,a-a_1)\mu_1 + a_1 \mu_1' + \min(b-b_1,k_2-a_1)\mu_2 \\ & + b_1 \mu_2' + \min(c-c_1,k_3-b_1)\mu_3 + c_1 \mu_3' + \min(d,k_4-c_1)\mu_4 + \max(0,a-(k_1 + a_1))\theta_1 \\ &+ \max(0,b-(k_2-a_1+b_1))\theta_2 + \max(0,c-(k_3-b_1+c_1))\theta_3\\ &+\max(0,d-(k_4-c_1))\theta_4+s
    \end{split}
\end{equation}
and
$\Tilde{C}_{a,a_1,b,b_1,c,c_1,d}(\mathbf{\Theta},s) = 0$ for $(a,a_1,b,b_1,c,c_1,d) \notin S$. 

\newpage
\section{Expected time spent in each state during the time interval \texorpdfstring{${\bm{(0,T)}}$}{(0,T)} -- four levels}\label{Ts1}

The Laplace transform of the expected time spent in each state in the case four levels of calls, ${\Tilde{\bm{v}}}_{a,a_1,b,b_1,c,c_1,d}(\bm{\Theta},s)$ is given by
\begin{equation}
{\Tilde{\bm{v}}}_{a,a_1,b,b_1,c,c_1,d}(\bm{\Theta},s) = \dfrac{\bm{\xi}^s_{a,a_1,b,b_1,c,c_1,d}(\bm{\Theta},s)}{\psi^s_{a,a_1,b,b_1,c,c_1,d}(s)},
\end{equation}
where
\begin{equation}
    \begin{split}
        \xi^s_{a,a_1,b,b_1,c,c_1,d}(\bm{\Theta},s) = &\bm{e_{a,a_1,b,b_1,c,c_1,d}}/s + (\lambda_1 {\Tilde{\bm{v}}}_{a+1,a_1+\mathbbm{1}(a-a_1 \geq k_1, b - b_1 < k_2 - \Theta_2 - a_1),b,b_1,c,c_1,d}(\bm{\Theta},s) \\&+ \lambda_2 {\Tilde{\bm{v}}}_{a,a_1,b+1,b_1+\mathbbm{1}(b-b_1 \geq k_2, c < k_3 - \Theta_3 - b_1),c,c_1,d}(\bm{\Theta},s)\\&+ \lambda_3 {\Tilde{\bm{v}}}_{a,a_1,b,b_1,c+1,c_1+\mathbbm{1}(c-c_1 \geq k_3, d < k_4 - \Theta_4 - c_1),d}(\bm{\Theta},s)\\ & + \lambda_4 {\Tilde{\bm{v}}}_{a,a_1,b,b_1,c,c_1,d+1}(\bm{\Theta},s)) \mathbbm{1}_{a+b+c+d < \ell}+\min(k_1,a-a_1)\mu_1 {\Tilde{\bm{v}}}_{a-1,a_1,b,b_1,c,c_1,d}(\bm{\Theta},s) \\ &+ a_1 \mu_1' {\Tilde{\bm{v}}}_{a-1,a_1 - (1-\mathbbm{1}(a-a_1 > k_1, b - b_1 = k_2 - \Theta_2 - a_1),b,b_1,c,c_1,d}(\bm{\Theta},s))\\ &+\min(b-b_1,k_2-a_1)\mu_2 {\Tilde{\bm{v}}}_{a,a_1+\mathbbm{1}(a-a_1 > k_1, b - b_1 = k_2 - \Theta_2 - a_1),b-1,b_1,c,c_1,d}(\bm{\Theta},s)\\ & + b_1 \mu_2' {\Tilde{\bm{v}}}_{a,a_1,b-1,b_1 - (1-\mathbbm{1}(b-b_1 > k_2 - a_1, c = k_3 - \Theta_3 - b_1),c,c_1,d}(\bm{\Theta},s)) \\&+\min(c-c_1,k_3-b_1)\mu_3 {\Tilde{\bm{v}}}_{a,a_1,b,b_1+\mathbbm{1}(b-b_1 > k_2, c - c_1 = k_3 - \Theta_3 - b_1),c-1,c_1,d}(\bm{\Theta},s)\\ & + c_1 \mu_3' {\Tilde{\bm{v}}}_{a,a_1,b,b_1,c-1,c_1 - (1-\mathbbm{1}(c-c_1 > k_3 - b_1, d = k_4 - \Theta_4 - c_1),d}(\bm{\Theta},s)) \\ &+ \min(d,k_4-c_1)\mu_4 {\Tilde{\bm{v}}}_{a,a_1,b,b_1,c,c_1+\mathbbm{1}(c-c_1 > k_3 - b_1, d = k_4 - \Theta_4 - c_1),d-1}(\bm{\Theta},s) \\ &+\max(0,a-(k_1 + a_1))\theta_1 {\Tilde{\bm{v}}}_{a-1,a_1,b,b_1,c,c_1,d}(\bm{\Theta},s)\\ &+\max(0,b-(k_2-a_1+b_1))\theta_2 {\Tilde{\bm{v}}}_{a,a_1,b-1,b_1,c,c_1,d}(\bm{\Theta},s)\\ &+\max(0,c-(k_3-b_1+c_1))\theta_3 {\Tilde{\bm{v}}}_{a,a_1,b,b_1,c-1,c_1,d}(\bm{\Theta},s) \\ &+\max(0,d-(k_4-c_1))\theta_4 {\Tilde{\bm{v}}}_{a,a_1,b,b_1,c,c_1,d-1}(\bm{\Theta},s),
    \end{split}
\end{equation}
\begin{equation}
    \begin{split}
        \psi^s_{a,a_1,b,b_1,c,c_1,d}(s) = &(\lambda_1+\lambda_2+\lambda_3+\lambda_4)\mathbbm{1}_{a+b+c+d < \ell} + \min(k_1,a-a_1)\mu_1 + a_1 \mu_1' + \min(b-b_1,k_2-a_1)\mu_2 \\ & + b_1 \mu_2' + \min(c-c_1,k_3-b_1)\mu_3 + c_1 \mu_3' + \min(d,k_4-c_1)\mu_4 + \max(0,a-(k_1 + a_1))\theta_1 \\ &+ \max(0,b-(k_2-a_1+b_1))\theta_2 + \max(0,c-(k_3-b_1+c_1))\theta_3\\ &+\max(0,d-(k_4-c_1))\theta_4+s,
    \end{split}
\end{equation}
and
${\Tilde{\bm{v}}}_{a,a_1,b,b_1,c,c_1,d}(\bm{\Theta},s) = 0$ for $(a,a_1,b,b_1,c,c_1,d) \notin S$.

\newpage
\section{Waiting time distributions -- four levels}

\subsection{Level 2 callers}\label{l2ap}

The Laplace transform of ${\Xi}_{a_1,b,b_1,c,c_1,d}^2(\bm{\Theta},y)$ given by $\tilde{\Xi}_{a_1,b,b_1,c,c_1,d}^2(\bm{\Theta},s)$ is equal to $0$ for $a_1 + b + c + d = \ell$ and otherwise,
\begin{equation}\tilde{\Xi}_{a_1,b,b_1,c,c_1,d}^2(\bm{\Theta},s) = \begin{cases}
    \dfrac{1}{s}, & k_2 - a_1 > b - b_1 \text{ or } k_3-\Theta_3 > c + b_1 \\[3ex]
    \dfrac{\xi^2_{a_1,b,b_1,c,c_1,d}(\bm{\Theta},s)}{\psi^2_{a_1,b,b_1,c,c_1,d}(s)}, & k_2 - a_1 \leq b - b_1  \text{ and } k_3-\Theta_3 \leq c + b_1
\end{cases}\end{equation}
where
\begin{equation}
    \begin{split}
        \xi^2_{a_1,b,b_1,c,c_1,d}(\bm{\Theta},s) = & (\lambda_3 \tilde{\Xi}_{a_1,b,b_1,c+1,c_1,d}^2 + \lambda_4 \tilde{\Xi}_{a_1,b,b_1,c,c_1,d+1}^2)\mathbbm{1}_{a_1+b+c+d < \ell-1} + \\ & a_1 \mu_1' \tilde{\Xi}_{a_1-1,b,b_1,c,c_1,d}^2(\bm{\Theta},s) + (k_2 -a_1)\mu_2 \tilde{\Xi}_{a_1,b-1,b_1,c,c_1,d}^2(\bm{\Theta},s) \\ & + b_1\mu_2'\tilde{\Xi}_{a_1,b-1,b_1-(1-\mathbbm{1}_{b+1-b_1 > k_2-a_1, c = k_3 - b_1 - \Theta_3}),c,c_1,d}^2(\bm{\Theta},s) \\ & + \min(k_3-b_1,c-c_1)\mu_3 \tilde{\Xi}_{a_1,b,b_1 + \mathbbm{1}_{b+1-b_1 > k_4-c_1, c-c_1 = k_3-b_1-\Theta_3},c-1,c_1,d}^2(\bm{\Theta},s) \\ &+ c_1\mu_3'\tilde{\Xi}_{a_1,b,b_1,c-1,c_1-(1-\mathbbm{1}_{c-c_1 > k_3-b_1, d = k_4 - c_1 - \Theta_4}),d}^2(\bm{\Theta},s) \\ & + \min(k_4-c_1,d)\mu_4 \tilde{\Xi}_{a_1,b,b_1,c,c_1 + \mathbbm{1}_{c-c_1 > k_3-b_1, d = k_4-c_1-\Theta_4},d-1}^2(\bm{\Theta},s) \\ & +(b-(k_2-a_1+b_1))\theta_2 \tilde{\Xi}_{a_1,b-1,b_1,c,c_1,d}^2(\bm{\Theta},s) + \\ & \max(0,c-(k_3-b_1+c_1))\theta_3\tilde{\Xi}_{a_1,b,b_1,c-1,c_1,d}^2(\bm{\Theta},s) \\ & \max(0,d-(g_4-c_1))\theta_4 \tilde{\Xi}_{a_1,b,b_1,c,c_1,d-1}^2(\bm{\Theta},s) 
    \end{split}
\end{equation}
and
\begin{equation}
    \begin{split}
        \psi^2_{a_1,b,b_1,c,c_1,d}(s) = &(\lambda_3 + \lambda_4) \mathbbm{1}_{a_1+b+c+d < \ell-1} + a_1 \mu_1' + (k_2 -a_1)\mu_2 + b_1\mu_2' + c_1\mu_3' \\ & + \min(c-c_1, k_3-b_1)\mu_3 + \min(d,k_4-c_1)\mu_4 +(b+1-(k_2-a_1+b_1))\theta_2 \\ &+ \max(0,c-(k_3-b_1+c_1))\theta_3 + \max(0,d-(k_4-c_1))\theta_4 + s.
    \end{split}
\end{equation}

\subsection{Level 1 callers}\label{l1ap}

The Laplace transform of ${\Xi}_{a,a_1,b,b_1,c,c_1,d}^1(\bm{\Theta},y)$ denoted by $\tilde{\Xi}_{a,a_1,b,b_1,c,c_1,d}^1(\bm{\Theta},s)$ is given by
\begin{equation}\tilde{\Xi}_{a,a_1,b,b_1,c,c_1,d}^1(\bm{\Theta},s) = \begin{cases}
    \dfrac{1}{s}, & k_1 > a - a_1 \text{ or } k_2-\Theta_2 > b + a_1 \\[3ex]
    \dfrac{\xi^1_{a,a_1,b,b_1,c,c_1,d}(\bm{\Theta},s)}{\psi^1_{a,a_1,b,b_1,c,c_1,d}(s)}, & k_1 \leq a - a_1  \text{ and } k_2-\Theta_2 \leq b + a_1
\end{cases}\end{equation}
where
\begin{equation}
    \begin{split}
        \xi^1_{a,a_1,b,b_1,c,c_1,d}(\bm{\Theta},s) = & (\lambda_2 \Tilde{\Xi}_{a,a_1,b+1,b_1+\mathbbm{1}(b-b_1 \geq k_2, c < k_3 - \Theta_3 - b_1),c,c_1,d}^1(\bm{\Theta},s)\\&+ \lambda_3 \Tilde{\Xi}_{a,a_1,b,b_1,c+1,c_1+\mathbbm{1}(c-c_1 \geq k_3, d < k_4 - \Theta_4 - c_1),d}^1(\bm{\Theta},s)\\ & + \lambda_4 \Tilde{\Xi}_{a,a_1,b,b_1,c,c_1,d+1}^1(\bm{\Theta},s)) \mathbbm{1}_{a+b+c+d < \ell-1}+k_1\mu_1 \Tilde{\Xi}_{a-1,a_1,b,b_1,c,c_1,d}^1(\bm{\Theta},s) \\ &+ a_1 \mu_1' \Tilde{\Xi}_{a-1,a_1 - (1-\mathbbm{1}(a-a_1 > k_1, b - b_1 = k_2 - \Theta_2 - a_1),b,b_1,c,c_1,d}^1(\bm{\Theta},s))\\ &+\min(b-b_1,k_2-a_1)\mu_2 \Tilde{\Xi}_{a,a_1+\mathbbm{1}(a-a_1 > k_1, b - b_1 = k_2 - \Theta_2 - a_1),b-1,b_1,c,c_1,d}^1(\bm{\Theta},s)\\ & + b_1 \mu_2' \Tilde{\Xi}_{a,a_1,b-1,b_1 - (1-\mathbbm{1}(b-b_1 > k_2 - a_1, c = k_3 - \Theta_3 - b_1),c,c_1,d}^1(\bm{\Theta},s)) \\&+\min(c-c_1,k_3-b_1)\mu_3 \Tilde{\Xi}_{a,a_1,b,b_1+\mathbbm{1}(b-b_1 > k_2, c - c_1 = k_3 - \Theta_3 - b_1),c-1,c_1,d}^1(\bm{\Theta},s)\\ & + c_1 \mu_3' \Tilde{\Xi}_{a,a_1,b,b_1,c-1,c_1 - (1-\mathbbm{1}(c-c_1 > k_3 - b_1, d = k_4 - \Theta_4 - c_1),d}^1(\bm{\Theta},s)) \\ &+ \min(d,k_4-c_1)\mu_4 \Tilde{\Xi}_{a,a_1,b,b_1,c,c_1+\mathbbm{1}(c-c_1 > k_3 - b_1, d = k_4 - \Theta_4 - c_1),d-1}^1(\bm{\Theta},s) \\ &+(a-(k_1 + a_1))\theta_1 \Tilde{\Xi}_{a-1,a_1,b,b_1,c,c_1,d}^1(\bm{\Theta},s)\\ &+\max(0,b-(k_2-a_1+b_1))\theta_2 \Tilde{\Xi}_{a,a_1,b-1,b_1,c,c_1,d}^1(\bm{\Theta},s)\\ &+\max(0,c-(k_3-b_1+c_1))\theta_3 \Tilde{\Xi}_{a,a_1,b,b_1,c-1,c_1,d}^1(\bm{\Theta},s) \\ &+\max(0,d-(k_4-c_1))\theta_4 \Tilde{\Xi}_{a,a_1,b,b_1,c,c_1,d-1}^1(\bm{\Theta},s)
    \end{split}
\end{equation}
and
\begin{equation}
    \begin{split}
        \psi^1_{a,a_1,b,b_1,c,c_1,d}(s) = &(\lambda_2+\lambda_3+\lambda_4)\mathbbm{1}_{a+b+c+d < \ell-1} + k_1\mu_1 + a_1 \mu_1' + b_1 \mu_2' \\ & + \min(b-b_1,k_2-a_1)\mu_2 + \min(c-c_1,k_3-b_1)\mu_3 + c_1 \mu_3' \\ & + \min(d,k_4-c_1)\mu_4 + (a+1-(k_1 + a_1))\theta_1 \\ & + \max(0,b-(k_2-a_1+b_1))\theta_2 + \max(0,c-(k_3-b_1+c_1))\theta_3 \\ & + \max(0,d-(k_4-c_1))\theta_4+s,
    \end{split}
\end{equation}
for $a+b+c+d < \ell$, and otherwise $\tilde{\Xi}_{a,a_1,b,b_1,c,c_1,d}^1(\bm{\Theta},s) = 0$ for $a+b+c+d = \ell$.

\end{document}